\documentclass[a4paper,12pt]{article}
\usepackage[latin1]{inputenc}
\usepackage[francais,english]{babel}    
\usepackage{amssymb}
\usepackage{amsmath}
\usepackage{amsthm}
\usepackage{multicol}
\usepackage{graphicx}
\usepackage{color}
\usepackage[T1]{fontenc}
\usepackage{yfonts}
\usepackage{placeins}
\usepackage{array}
\usepackage{here}

\newtheorem{thm}{Théorème}[section]
\newcommand{\bthm}{\begin{thm}}
\newcommand{\ethm}{\end{thm}}

\newtheorem*{thn}{Théorème}
\newcommand{\bthn}{\begin{thn}}
\newcommand{\ethn}{\end{thn}}

\newtheorem{pro}[thm]{Proposition}
\newcommand{\bpro}{\begin{pro}}
\newcommand{\epro}{\end{pro}}

\newtheorem{cor}[thm]{Corollaire}
\newcommand{\bcor}{\begin{cor}}
\newcommand{\ecor}{\end{cor}}

\newtheorem{lem}[thm]{Lemme}
\newcommand{\blem}{\begin{lem}}
\newcommand{\elem}{\end{lem}}

\newtheorem*{rmq}{Remarque}
\newcommand{\brq}{\begin{rmq} \upshape}
\newcommand{\erq}{\end{rmq}}

\newtheorem*{exe}{Exemple}
\newcommand{\bexe}{\begin{exe} \upshape}
\newcommand{\eexe}{\end{exe}}

\newtheorem*{pre}{Démonstration}
\newcommand{\bp}{\begin{pre} \upshape}
\newcommand{\ep}{\hfill \qed \end{pre}}

\newtheorem*{pthm}{Preuve du théorème~\ref{thm:unicite_ecriture}}
\newcommand{\bpthm}{\begin{pthm}}
\newcommand{\epthm}{\end{pthm}}
\newcommand{\beq}{\begin{eqnarray*}}
\newcommand{\eeq}{\end{eqnarray*}}

\newcommand{\beqn}{\begin{equation}}
\newcommand{\eeqn}{\end{equation}}

\newcommand{\ben}{\begin{enumerate}}
\newcommand{\een}{\end{enumerate}}
\newcommand{\bit}{\begin{itemize} \renewcommand{\labelitemi}{$\bullet$} }
\newcommand{\eit}{\end{itemize}}

\newcommand{\bfg}{
\begin{figure}[H]
\begin{center}}
\newcommand{\efg}{
\end{center}
\end{figure}
\FloatBarrier}

\newcolumntype{M}[1]{>{\raggedright}m{#1}}

\newcommand{\df}{\emph}

\newcommand{\R}{\mathbb{R}}

\newcommand{\N}{\mathbb{N}}
\newcommand{\Z}{\mathbb{Z}}
\newcommand{\C}{\mathbb{C}}

\renewcommand{\SS}{\mathbb{S}}
\renewcommand{\H}{\mathbb{H}}

\newcommand{\bs}{\symbol{92}}
\newcommand{\Ch}{\mathcal{S}}

\newcommand{\ov}{\overline}

\renewcommand{\Im}{\operatorname{Im}}

\newcommand{\Ker}{\operatorname{Ker}}
\newcommand{\Vect}{\operatorname{Vect}}

\renewcommand{\t}{ ^t\!}
\newcommand{\tr}{\operatorname{tr}}

\newcommand{\Diag}{\operatorname{Diag}}
\renewcommand{\dim}{\operatorname{dim}}
\newcommand{\Card}{\operatorname{Card}}
\newcommand{\Isom}{\operatorname{Isom}}
\newcommand{\Id}{\operatorname{Id}}
\newcommand{\Lie}{\operatorname{Lie}}
\newcommand{\CAT}{\operatorname{CAT}}
\newcommand{\Sp}{\operatorname{Sp}}

\newcommand{\rang}{\operatorname{rang}}

\newcommand{\eps}{\varepsilon}

\newcommand{\ra}{\rightarrow}
\newcommand{\ral}[1]{\underset{#1}{\longrightarrow}}

\newcommand{\liml}{\lim\limits}

\renewcommand{\geq}{\geqslant}
\renewcommand{\leq}{\leqslant}

\newcommand{\gothique}{\mathfrak}
\renewcommand{\gg}{\gothique{g}}
\newcommand{\kk}{\textgoth{k}}
\newcommand{\pp}{\gothique{p}}
\renewcommand{\aa}{\gothique{a}}
\newcommand{\ddd}{\gothique{d}}

\newcommand{\zz}{\gothique{z}}
\newcommand{\mm}{\gothique{m}}
\newcommand{\nn}{\gothique{n}}

\newcommand{\sll}{\gothique{sl}}
\newcommand{\so}{\gothique{so}}
\newcommand{\gl}{\gothique{gl}}

\renewcommand{\exp}{\operatorname{exp}}
\renewcommand{\log}{\operatorname{log}}

\newcommand{\ad}{\operatorname{ad}}
\newcommand{\Ad}{\operatorname{Ad}}

\newcommand{\GL}{\operatorname{GL}}
\newcommand{\SL}{\operatorname{SL}}
\newcommand{\SO}{\operatorname{SO}}
\renewcommand{\O}{\operatorname{O}}

\newcommand{\PSL}{\operatorname{PSL}}

\newcommand{\Sym}{\operatorname{Sym}}

\makeatletter
\def\Ddots{\mathinner{\mkern1mu\raise\p@
\vbox{\kern7\p@\hbox{.}}\mkern2mu
\raise4\p@\hbox{.}\mkern2mu\raise7\p@\hbox{.}\mkern1mu}}
\makeatother

\makeatletter
\def\maketitles{%
  \null
  \thispagestyle{empty}%
  \vfill
  \begin{center}\leavevmode
    \normalfont
    {\LARGE \@title\par}%
    \vskip 1.2cm
    {\large \@author\par}%
    \vskip 1.2cm
    {\large \@subtitle\par}%
    \vskip 0.8cm
    {\large \@date\par}%
  \end{center}%
  \vfill
  \null
  \cleardoublepage
  }

\def\date#1{\def\@date{#1}}
\def\author#1{\def\@author{#1}}
\def\title#1{\def\@title{#1}}
\def\subtitle#1{\def\@subtitle{#1}}
\makeatother

\title{Compactification de Chabauty des espaces sym\'etriques de type non compact}                                    
\author{Thomas Haettel}                 
\date{10 Mai 2010}  

\begin{document}

\maketitle

\begin{abstract}
The space of closed subgroups of a locally compact topological group is endowed with a natural topology, called the Chabauty topology. Let $X$ be a symmetric space of noncompact type, and $G$ be its group of isometries. The space $X$ identifies with the subspace of maximal compact subgroups of $G$ : taking the closure gives rise to the Chabauty compactification of the symmetric space $X$. Using simpler arguments than those present in~(\cite{guivarch} Guivarc'h, Y., L.~Ji and J.C.~Taylor, ``Compactifications of symmetric spaces", Progr.~Math. {\bf 156}, 1998.), we describe the subgroups that appear in the boundary of the compactification, and classify the maximal distal and maximal amenable subgroups of $G$. We also provide a straightforward identification between the Chabauty compactification and the polyhedral compactification.
\end{abstract}

\selectlanguage{french}

\section*{Introduction}

\addcontentsline{toc}{section}{Introduction}

Si $G$ est un groupe topologique localement compact, munissons l'ensemble $\Ch(G)$ des sous-groupes fermés de $G$ de la topologie de Chabauty (voir par exemple~\cite{chabauty}, \cite[Chap.~VIII, $\S 5$]{bourbaki}, \cite{harpe_chabauty} et \cite{cdp}, ainsi que la partie~\ref{sec:chabauty}). Cette topologie fait de $\Ch(G)$ un espace compact, qui est métrisable pour la distance de Hausdorff pointée si $G$ est métrisable. L'espace $\Ch(G)$ est en général difficile à expliciter : si l'espace $\Ch(\C)$ est homéomorphe à la sphère $\SS^4$ (voir~\cite{pourezza}), pour d'autres petits groupes $G$, la topologie de l'espace $\Ch(G)$ est bien plus complexe (voir par exemple~\cite{heisenberg}, \cite{kloeckner} et \cite{RxZ}).

\bigskip

Si $X$ est un espace topologique localement compact et $G$ un groupe topologique localement compact agissant continûment sur $X$, considérons l'application de $X$ dans $\Ch(G)$ qui à un point de $X$ associe son stabilisateur. Si $\widehat{X} = X \cup \{\infty\}$ désigne le compactifié d'Alexandrov de $X$, identifions $X$ avec son image par le plongement diagonal dans $\Ch(G) \times \widehat{X}$. L'adhérence $\ov{X}^\Ch$ de (l'image de) $X$ est appelée la \df{compactification de Chabauty} de $X$ définie par cette action de $G$. Les groupes apparaissant dans la frontière de $X$ dans $\ov{X}^\Ch$ seront appelés les \df{groupes limites}.

\bigskip

Nous nous intéressons dans cet article au cas où $X$ est un espace symétrique de type non compact, comme par exemple les espaces hyperboliques réels ou les variétés riemanniennes quotients $\SL(n,\R) / \SO(n,\R)$ (voir par exemple~\cite{helgason} et la partie~\ref{sec:espaces} pour des rappels). Les compactifications de ces espaces ont été très étudiées (voir par exemple les travaux de Mostow~\cite{mostow}, Satake~\cite{satake}, Borel et Ji~\cite{borel_ji}...). Elles sont nombreuses, en fonction des besoins nécessaires : compactifications géodésique, de Satake, de Furstenberg, polyédrale, de Martin, de Karpelevic, etc. (voir par exemple~\cite{guivarch}). La géométrie des espaces symétriques étant intimement liée à la structure de leurs groupes d'isométries, il était naturel, comme dans~\cite{guivarch}, d'étudier leurs compactifications de Chabauty.

\bigskip

Soit $G$ un groupe de Lie semi-simple connexe, de centre fini et sans facteur compact. Nous renvoyons par exemple à~\cite{helgason} et à la partie~\ref{sec:espaces} pour des rappels, en particulier sur la terminologie qui suit. Soit $K$ un sous-groupe compact maximal de $G$, une décomposition d'Iwasawa $G=KAN$ de $G$ et $X=G/K$, muni d'une métrique riemannienne $G$-invariante, l'espace symétrique de type non compact associé. L'algèbre de Lie $\aa$ de $A$ est une sous-algèbre de Lie abélienne maximale de l'algèbre de Lie $\gg$ de $G$. Elle définit un système de racines $\Sigma$, qui est l'ensemble des formes linéaires non nulles $\alpha$ sur $\aa$ dont l'espace de poids $\gg_\alpha=\{Y \in \gg : \forall H \in \aa, \ad H(Y) = \alpha(H)Y\}$ est non nul. Choisissons une base $\Delta$ du système de racines défini par $\aa$, et notons $\ov{\aa^+}$ la chambre de Weyl fermée associée. Notons de plus $M=Z_K(\aa)$ le centralisateur de $\aa$ dans $K$.

Pour toute partie propre $I$ de $\Delta$, notons $\aa_I = \bigcap_{\alpha \in I} \Ker \alpha$ et $\aa^I$ l'orthogonal de $\aa_I$ dans $\aa$ pour la forme de Killing de $G$. Notons $A^I$ le sous-groupe de Lie connexe de $G$ d'algèbre de Lie $\aa^I$, et $\ov{A^{I,+}} = A^I \cap \exp \ov{\aa^+}$.

Définissons alors le sous-groupe de Lie $G^I = D(Z(\aa_I))_0$, composante neutre du groupe dérivé du centralisateur de $\aa_I$ dans $G$, et $K^I = G^I \cap K$. Soit $\Sigma_I^+$ l'ensemble des racines positives qui ne sont pas combinaison linéaire des racines de $I$, et soit $\nn_I$ la sous-algèbre de Lie somme directe des espaces de poids de $\Sigma_I^+$. Notons $N_I$ le sous-groupe de Lie connexe de $G$ d'algèbre de Lie $\nn_I$. Pour tous $a \in A$ et $k \in K$, notons
$$ D^I_{a,k} = kaK^IMN_Ia^{-1}k^{-1}.$$

\bigskip

Le but principal de cet article est de donner une nouvelle preuve, plus courte et plus directe, du théorème de~\cite{guivarch} donnant une description explicite des groupes limites de la compactification de Chabauty de $X$.

\bthn Soit $D \in \ov{X}^\Ch \bs X$ un groupe limite. Alors il existe une partie propre $I$ de $\Delta$, $a \in \ov{A^{I,+}}$ et $k \in K$ tels que $D=D^I_{a,k}$. 

De plus, cette écriture est unique au sens suivant : pour $I_1,I_2$ deux parties propres de $\Delta$, $a_1 \in \ov{A^{I_1,+}}$, $a_2 \in \ov{A^{I_2,+}}$ et $k_1,k_2 \in K$, nous avons $D^{I_1}_{a_1,k_1}=D^{I_2}_{a_2,k_2}$ si et seulement si
$$ I_1=I_2=I, a_1 = a_2 =a \mbox{ et } k_2^{-1}k_1 \in (K^I \cap aK^Ia^{-1})M. $$
\ethn

Par exemple, considérons le groupe $G=\SL(2,\R)^n$, dont l'espace symétrique associé est le $n$-polydisque $\H_2^n$. Pour $n=1$, la compactification de Chabauty $\ov{\H_2}^\Ch$ coïncide avec la compactification géométrique $\ov{\H_2}^g$ du disque obtenue en ajoutant un cercle à l'infini. Pour $n \geq 2$, on voit que la compactification de Chabauty est le produit des $n$ compactifications de Chabauty, c'est-à-dire le produit de $n$ disques fermés : $\ov{\H_2^n}^\Ch = \left(\ov{\H_2}^\Ch\right)^n$. En effet, les groupes limites de $G$ sont les produits de $n$ groupes limites de $\SL(2,\R)$. En revanche, la compactification géométrique du $n$-polydisque est obtenue en ajoutant une $(2n-1)$-sphère à l'infini.

\bigskip

Nous retrouvons ensuite les résultats (voir~\cite{moore_amenable} et \cite{guivarch_eigen}) de classification des sous-groupes distaux maximaux et moyennables maximaux de $G$ : la compactification $\ov{X}^\Ch$ est l'ensemble des sous-groupes distaux maximaux de $G$, et les sous-groupes moyennables maximaux de $G$ sont les normalisateurs de ces sous-groupes.

\bigskip
La première partie de cet article est constituée de rappels bien connus sur les décom\-positions des groupes de Lie semi-simples, et pose les notations. La deuxième partie présente brièvement la topologie de Chabauty. Le c\oe ur de l'article est la partie $3$, en particulier le théorème~\ref{thm:existence_ecriture}~: tandis que la preuve originelle de~\cite{guivarch} est fondée sur l'étude de mesures et l'utilisation des frontières de Furstenberg, la preuve présentée ici n'utilise que des arguments élémentaires de groupes de Lie. Celle-ci semble adaptée à l'étude de la compactification de Chabauty, qui présente l'avantage par rapport aux autres compactifications d'avoir une définition élémentaire. Enfin, la dernière partie décrit l'homéomorphisme entre les compactifications de Chabauty et polyédrale. L'existence de cet homéomorphisme est bien connue, mais cette correspondance ne semblait pas exister dans la littérature. Par ailleurs, la continuité de l'action de $G$ sur la compactification polyédrale est montrée directement, alors que dans~\cite{guivarch} elle est montrée via la compactification de Martin.

\bigskip

Je tiens à remercier chaleureusement Frédéric Paulin pour sa disponibilité et ses nombreux conseils, ainsi que le rapporteur pour ses suggestions.

\section{Décompositions des espaces symétriques de type non compact}

\label{sec:espaces}

\bigskip
\textbf{Préliminaires}
\bigskip

Toutes les algèbres de Lie et tous les groupes de Lie considérés sont, sauf mention contraire, réels. Concernant les espaces symétriques, une référence assez complète est~\cite{helgason}.

\bigskip

Une variété riemannienne connexe $X$ est appelée un \df{espace (globalement) symétrique} si tout point $x$ de $X$ est un point fixe isolé d'une isométrie involutive de $X$. Il est dit de plus \df{de type non compact} si sa courbure sectionnelle est partout négative ou nulle, et si le revêtement universel de $X$ n'admet pas de facteur de De Rham euclidien non trivial.

\bexe Les espaces hyperboliques réels sont des espaces symétriques de type non compact. \eexe

Rappelons que si $X$ est un espace symétrique de type non compact, alors la composante neutre $G=\Isom_0(X)$ du groupe des isométries de $X$ est un groupe de Lie connexe semi-simple, de centre trivial et sans facteur compact (voir par exemple~\cite[Exercise~5, p.~250]{helgason}). De plus, si l'on choisit un point base $x$ de $X$, alors l'évaluation en $x$ induit un difféomorphisme de $G/G_x$ sur $X$, où l'on note $G_x$ le stabilisateur de $x$ dans $G$ (voir par exemple~\cite[Theorem~3.3, p.~208]{helgason}). Le groupe $G_x$ est un sous-groupe compact maximal de $G$, et par ailleurs les sous-groupes compacts maximaux de $G$ sont les stabilisateurs des points de $X$, et sont deux à deux conjugués (voir par exemple~\cite[Theorem~3, p.~259]{onishchik_vinberg}).

Réciproquement, si $G$ est un groupe de Lie connexe semi-simple, de centre fini et sans facteur compact et si $K$ est un sous-groupe compact maximal de $G$, alors toute métrique riemannienne $G$-invariante sur l'espace $G/K$ en fait un espace symétrique de type non compact, et le groupe $G$ se surjecte canoniquement sur la composante neutre $\Isom_0(G/K)$ du groupe des isométries de $G/K$, avec pour noyau le centre fini de $G$ (voir par exemple~\cite[Proposition~3.4, p.~209 et Theorem~3.1, p.~241]{helgason}).

\bexe Pour tout $n \geq 2$, l'espace $\SL(n,\R)/\SO(n,\R)$, muni d'une métrique\-$\SL(n,\R)$-invariante, est un espace symétrique de type non compact. \eexe

Soit $G$ un groupe de Lie semi-simple de centre fini ayant un nombre fini de composantes connexes. Choisissons un sous-groupe compact maximal $K$ de $G$. Alors il existe une unique décomposition de Cartan $\gg=\kk \oplus \pp$ de $\gg$, où $\gg$ et $\kk$ désignent les algèbres de Lie de $G$ et $K$ respectivement, orthogonale pour la forme de Killing $B:(X,Y) \ra \tr(\ad X \circ \ad Y)$ de $\gg$. En outre, rappelons que si $\theta$ désigne l'involution de Cartan de $G$ associée, valant l'identité sur $\kk$ et moins l'identité sur $\pp$, alors la forme bilinéaire $B_\theta : (Y,Z) \mapsto -B(Y,\theta(Z))$ est définie positive sur $\gg$. En tant qu'endomorphismes de l'espace vectoriel $\gg$ muni du produit scalaire $B_\theta$, les éléments de $\ad \kk$ sont antisymétriques et ceux de $\ad \pp$ sont symétriques.

\bexe Pour le sous-groupe compact maximal $\SO(n,\R)$ de $\SL(n,\R)$, la décom\-position de Cartan est $\sll(n,\R) = \so(n,\R) \oplus \Sym_0(n,\R)$, où $\so(n,\R)$ est la sous-algèbre de Lie des matrices antisymétriques, et où $\Sym_0(n,\R)$ est l'espace vectoriel des matrices symétriques de trace nulle. \eexe

L'application $\pp \times K \rightarrow G$ définie par $(Y,k) \mapsto (\exp Y)k$ est un difféomorphisme, appelé décomposition polaire de $G$ (voir par exemple~\cite[Theorem~2, p.~256]{onishchik_vinberg}).

\bigskip

Choisissons $\aa$ une sous-algèbre de Lie de $\gg$ abélienne diagonalisable sur $\R$ maximale  (c'est-à-dire que les éléments de $\ad \aa$ sont diagonalisables sur $\R$). Elle sont toutes $\Ad G$-conjuguées entre elles, et on peut choisir $\aa$ incluse dans $\pp$. La dimension de $\aa$ est appelée le rang (réel) de l'algèbre de Lie $\gg$. Si l'on note $\aa^*$ l'espace vectoriel dual de $\aa$, définissons le système de racines (restreint) associé à $\aa$ :  
$$ \Sigma = \{ \alpha \in \aa^* \bs \{0\}  \,:\, \exists Y \in \gg \bs \{0\}, \forall H \in \aa, \ad H (Y) = \alpha(H) Y \}. $$
Ce système fournit une décomposition de $\gg$ en espaces de racines, orthogonale pour le produit scalaire $B_\theta$ :
$$\gg = \gg_0 \oplus \bigoplus_{\alpha \in \Sigma} \gg_\alpha $$
où on a posé $\gg_\alpha = \{Y \in \gg \,:\, \forall H \in \aa, \ad H (Y) = \alpha(H) Y\}$, et où $\gg_0=\zz_{\gg}(\aa)$ est le centralisateur dans $\gg$ de $\aa$. De plus, si l'on pose $\mm = \zz_{\kk}(\aa)$ le centralisateur dans $\kk$ de $\aa$, alors $\gg_0 = \mm \oplus \aa$.

\bexe Pour $\gg=\sll(n,\R)$, $\kk=\so(n,\R)$ et $\pp=\Sym_0(n,\R)$, on peut choisir pour $\aa$ la sous-algèbre des matrices diagonales. Dans ce cas, le système de racines est
$$\Sigma = \{\alpha_{i,j} \,:\, H \mapsto H_{i,i} - H_{j,j}, \forall i,j \in [[1,n]], i \neq j\}$$
et, pour $i,j \in [[1,n]]$, avec $i \neq j$, nous avons $\gg_{\alpha_{i,j}} = \Vect (E_{i,j})$, où $(E_{i,j})_{1 \leq i,j \leq n}$ désigne la base canonique de $\gl (n,\R)$. Et $\gg_0 = \aa$ dans ce cas particulier. \eexe

Rappelons comment cette décomposition se comporte par rapport au crochet de Lie : soient $\alpha,\beta \in \Sigma \cup \{0\}$. Alors
\beq [\gg_\alpha,\gg_\beta] \left\{ \begin{array}{l} \subset \gg_{\alpha+\beta} \mbox{ si } \alpha+\beta \in \Sigma \cup \{0\} \\ = \{0\} \mbox{ sinon.} \end{array} \right. \eeq

On appelle les composantes connexes de $\aa \bs \cup_{\alpha \in \Sigma} \Ker \alpha$ les \df{chambres de Weyl} (vectorielles, ouvertes) de $\aa$. Choisissons une chambre de Weyl $\aa^+$, que l'on appellera chambre de Weyl \df{positive}. On dit alors qu'une racine $\alpha \in \Sigma$ est \df{positive} si $\alpha|_{\aa^+} \geq 0$ (resp. \df{négative} si $\alpha|_{\aa^+} \leq 0$). Notons $\Sigma^+$ (resp. $\Sigma^- = \Sigma \bs \Sigma^+$) l'ensemble des racines positives (resp. négatives). Soit $\Delta$ l'ensemble des racines $\alpha \in \Sigma$ positives, ne pouvant pas s'écrire $\alpha = \beta+\gamma$, avec $\beta$ et $\gamma$ des racines positives. Alors $\Delta$ est une \df{base} du système de racines $\Sigma$, c'est-à-dire que toutes les racines de $\Sigma$ s'expriment de manière unique comme combinaison linéaire à coefficients entiers, tous de même signe, des racines de la base $\Delta$. De plus, toute base s'obtient ainsi, et on a
$$ \aa^+ = \{X \in \aa \,:\, \forall \alpha \in \Delta, \alpha(X) >0 \}.$$

\bexe Pour $G=\SL(n,\R)$, on peut choisir pour chambre de Weyl positive $\aa^+ = \{H \in \aa \,:\, H_{1,1} > H_{2,2} > ... > H_{n,n} \}$. Ceci correspond au choix de la base $\Delta = \{\alpha_{i,i+1}, i \in [[1,n-1]] \}$. \eexe

Posons $M=Z_K(\aa)$ le centralisateur de $\aa$ dans $K$ pour l'action adjointe, il a pour algèbre de Lie $\mm = \zz_\kk(\aa)$ le centralisateur de $\aa$ dans $\kk$. Notons $\nn$ et $\nn^-$ les sous-algèbres de Lie nilpotentes de $\gg$ définies par $\nn = \oplus_{\alpha \in \Sigma^+} \gg_\alpha$ et $\nn^- = \oplus_{\alpha \in \Sigma^-} \gg_\alpha$. On peut alors écrire la décomposition $\gg = \gg_0 \oplus \nn \oplus \nn^-$. Notons $A$, $N$ et $N^-$ les uniques sous-groupes de Lie de $G$ connexes, d'algèbres de Lie respectives $\aa$, $\nn$ et $\nn^-$. Alors l'application $K \times A \times N \rightarrow G$ définie par $(k,a,n) \mapsto kan$ est un difféomorphisme, appelé décomposition d'Iwasawa (voir par exemple~\cite[Theorem~6, p.~275]{onishchik_vinberg}).

\bexe Pour $G=\SL(n,\R)$, $M$ est égal au sous-groupe des matrices diagonales ayant des $\pm 1$ sur la diagonale. Et $N$ est égal au sous-groupe des matrices triangulaires supérieures de coefficients diagonaux égaux à $1$. \eexe

Définissons $A^+= \exp (\aa^+)$ la chambre de Weyl (ouverte) positive de $A$, et notons $\ov{\aa^+}$ et $\ov{A^+}$ (les adhérences dans $\gg$ et $G$) les chambres de Weyl fermées. Alors tout élément $g$ de $G$  s'écrit $g=k_1 a k_2$, où $k_1,k_2 \in K$ et $a \in \ov{A^+}$ : c'est la décomposition de Cartan $K\ov{A^+}K$. De plus, l'élément $a=a(g) \in \ov{A^+}$ est uniquement déterminé par $g$, et est continu en $g$ (voir par exemple~\cite[Theorem~1.1, p.~402]{helgason}).

\bigskip
\textbf{Décompositions à l'aide d'une partie de la base}
\bigskip

Choisissons une partie $I$ de la base $\Delta$ du système de racines $\Sigma$. Définissons alors le sous-espace vectoriel $\aa_I = \{H \in \aa \,:\, \forall \beta \in I, \beta(H) = 0\} = \bigcap_{\beta \in I} \Ker \beta$ de $\aa$, et $\aa^I$ l'orthogonal de $\aa_I$ dans $\aa$ pour la forme de Killing, qui est définie positive sur $\aa$. Remarquons que $\dim \aa^I = \Card I$.

\bexe Pour $G=\SL(n,\R)$, en prenant $I = \{\alpha_{i,i+1},\alpha_{j,j+1}\}$, où $1 \leq i < j \leq n-1$, nous avons
$$\aa_I = \left\{ \left( \begin{array}{ccc} a_1 & 0 & 0 \\ 0 & \ddots & 0 \\ 0 & 0 & a_n \end{array} \right) \in \aa \,:\, a_i=a_{i+1} \mbox{ et } a_j=a_{j+1} \right\}.$$
Si $i+1<j$, alors
$$\aa^I = \left\{ \left( \begin{array}{ccc} a_1 & 0 & 0 \\ 0 & \ddots & 0 \\ 0 & 0 & a_n \end{array} \right) \in \aa \,:\, a_i=-a_{i+1}, a_j=-a_{j+1} \mbox{ et } \forall k \not\in \{i,i+1,j,j+1\}, a_k=0 \right\}.$$
Si $i+1=j$, alors
$$\aa^I = \left\{ \left( \begin{array}{ccc} a_1 & 0 & 0 \\ 0 & \ddots & 0 \\ 0 & 0 & a_n \end{array} \right) \in \aa \,:\, a_i+a_{i+1}+a_{i+2}=0 \mbox{ et } \forall k \not\in \{i,i+1,i+2\}, a_k=0 \right\}.$$
\eexe

Ce choix d'une partie $I$ définit également une partition de l'ensemble des racines $\Sigma$ : soit $\Sigma^I$ l'ensemble des racines qui sont combinaison linéaire (à coefficients entiers) d'éléments de $I$.  Nous avons également
$$\Sigma^I = \{\alpha \in \Sigma : \forall H \in \aa_I, \alpha(H)=0\}.$$
Définissons de plus $\Sigma_I = \Sigma \bs \Sigma^I$ le complémentaire de $\Sigma^I$, ainsi que $\Sigma^{I,+} = \Sigma^I \cap \Sigma^+$ et $\Sigma_I^+ = \Sigma_I \cap \Sigma^+$.

\bigskip

Ceci définit également une décomposition de $\nn$ en somme directe $\nn = \nn^I \oplus \nn_I$, où $\nn^I = \oplus_{\alpha \in \Sigma^{I,+}} \gg_\alpha$ et $\nn_I = \oplus_{\alpha \in \Sigma_I^+} \gg_\alpha$. De plus, $\nn^I$ et $\nn_I$ sont des sous-algèbres de Lie de $\nn$, et $[\nn^I,\nn_I] \subset \nn_I$. Notons $A^I$, $A_I$, $N^I$ et $N_I$ les uniques sous-groupes de Lie immergés de $G$ connexes  d'algèbres de Lie respectives $\aa^I$, $\aa_I$, $\nn^I$ et $\nn_I$. Ils sont en fait fermés. Le groupe $A$ est égal au produit $A = A^I \times A_I$, et le groupe $N$ est égal au produit semi-direct $N = N_I \rtimes N^I$. Le groupe $A$ normalise $N_I$ et $N^I$, et les groupes $N^I$ et $A_I$ commutent. L'application exponentielle de chacun des groupes $A$, $A^I$, $A_I$, $N$, $N^I$ et $N_I$ est surjective (voir par exemple~\cite{m2} pour des détails).

\bexe Pour $G=\SL(n,\R)$, en prenant comme ci-dessus $I = \{\alpha_{i,i+1},\alpha_{j,j+1}\}$, où $1 \leq i < j \leq n-1$ :
$$\nn^I = \left\{ \left( \begin{array}{ccc} T_i & 0 & 0 \\0 & T_{j-i} & 0 \\ 0 & 0 & T_{n-j} \end{array} \right) \in \aa \right\} $$
$$\mbox{ et } \nn_I = \left\{ \left( \begin{array}{ccc} 0 & * & * \\ 0 & 0 & * \\ 0 & 0 & 0 \end{array} \right) \in \aa \right\}$$
où les blocs sont de dimensions respectives $i$, $j-i$ et $n-j$, et où $T_i \in \gl (i,\R)$ désigne une matrice triangulaire supérieure stricte quelconque.
\eexe

On définit de manière analogue $\Sigma^{I,-}$, $\Sigma_I^-$, $N^{I,-}$ et $N_I^-$, et nous avons alors des propriétés semblables.

\bigskip

Posons $G^I = D(Z_G(\aa_I))_0$ la composante neutre du groupe dérivé du centralisateur dans $G$ de $\aa_I$, ce groupe a pour algèbre de Lie $\gg^I = \ddd(\zz(\aa_I))$ l'algèbre de Lie dérivée du centralisateur de $\aa_I$. Notons également $K^I = G^I \cap K$, qui a pour algèbre de Lie $\kk^I = \gg^I \cap \kk$.

D'après~\cite[\S~2.13, p.~18]{guivarch}, le groupe de Lie $G^I$ est connexe, semi-simple, de centre fini et sans facteur compact, et $K^I$ en est un sous-groupe compact maximal. On peut choisir pour sous-algèbre de Cartan de $\gg^I$ l'algèbre de Lie $\aa^I$, dont le système de racines associé est $\Sigma^I|_{\aa^I}$ et dont on peut choisir pour base $I|_{\aa^I}$. La décomposition en espaces de racines de $\gg^I$ est
$$ \gg^I = \gg^I_0 \oplus \bigoplus_{\alpha \in \Sigma^I} \gg_\alpha.$$
De plus, la décomposition d'Iwasawa associée est $G^I=K^IA^IN^I$.

\bigskip
\textbf{L'exemple de $\SO_0(p,p)/\SO(p) \times \SO(p)$}
\bigskip

Soit $G=\SO_0(p,p)$ la composante neutre du groupe orthogonal d'une forme quadratique (réelle) de signature $(p,p)$, où $p$ est un entier supérieur ou égal à $2$. C'est un groupe de Lie connexe semi-simple de centre fini, égal à $\{\Id\}$ si $p$ est impair et à $\{\Id,-\Id\}$ si $p$ est pair.

\bigskip

Choisissons pour matrice de la forme quadratique la matrice carrée de taille $2p$ et d'ordre $2$
$$ J = \left( \begin{array}{ccc}
0 & 0 & 1 \\ 0 & \Ddots & 0 \\ 1 & 0 & 0 \end{array} \right) .$$

Alors les éléments de $G$, vu comme sous-groupe de $\SL(2p,\R)$, sont exactement les matrices $Y$ telles que $\t YJY = J$. De plus, l'algèbre de Lie $\gg = \so(p,p)$ de $G$ s'écrit
$$\gg = \{Y \in \sll(2p,\R) : J\t YJ = -Y \}.$$
Donc l'algèbre de Lie $\gg$ est constituée des matrices antisymétriques par rapport à la deuxième diagonale. En raison de cette symétrie, pour tout entier $i \in [[1,2p]]$, nous noterons $\ov{i} = 2p+1-i$ son "symétrique".

\bigskip

Posons $\theta : \gg \ra \gg$ le morphisme d'algèbres de Lie défini par $Y \mapsto - \t Y$, et notons $\kk = \gg \cap \so(2p)$ l'ensemble des points fixes par $\theta$ et $\pp = \gg \cap \Sym_0(2p)$ l'ensemble des vecteurs propres pour la valeur propre $-1$ (où $\Sym_0(2p)$ désigne le sous-espace vectoriel de $\gl(2p,\R)$ des matrices symétriques). Alors $\gg = \kk \oplus \pp$ est une décomposition de Cartan de $\gg$. La forme de Killing sur $\gg$ est donnée par $B(Y,Z) = 2p\tr(YZ)$.

\bigskip

Notons $\aa$ le sous-espace vectoriel de $\gg$ constitué des matrices diagonales : c'est une sous-algèbre de Lie abélienne maximale de $\gg$. Une base de l'espace vectoriel $\aa$ est donnée par les matrices $H_i = E_{i,i} - E_{\ov{i},\ov{i}}$, pour $i \in [[1,p]]$ (où $(E_{i,j})_{i,j}$ désigne la base canonique de $\gl(2p,\R)$). Notons $(\beta_i)_{i \in [[1,p]]}$ la base duale de cette base. La forme linéaire $\beta_i$ est l'application qui à une matrice diagonale $H \in \aa$ associe $H_{i,i}-H_{\ov{i},\ov{i}}$.

\bigskip

Définissons les vecteurs suivants, pour $1 \leq i < j \leq p$.
\beq X_{\beta_i-\beta_j} &=& E_{i,j} - E_{\ov{j},\ov{i}} \\
 X_{-(\beta_i-\beta_j)} &=& E_{j,i} - E_{\ov{i},\ov{j}} \\
 X_{\beta_i+\beta_j} &=& E_{i,\ov{j}} - E_{j,\ov{i}} \\
 X_{-(\beta_i+\beta_j)} &=& E_{\ov{j},i} - E_{\ov{i},j} \eeq
Alors le vecteur $X_\alpha$ est de poids $\alpha$. Le système de racines $\Sigma$ de $\gg$ associé à la sous-algèbre abélienne maximale $\aa$ est
$$\Sigma = \{\beta_i-\beta_j : i,j \in [[1,p]], i \neq j\} \cup \{\beta_i+\beta_j : i,j \in [[1,p]], i \neq j\}. $$
L'espace de poids correspondant à une racine $\alpha \in \Sigma$ est $\gg_\alpha = \R X_\alpha$. De plus, le centralisateur de $\aa$ dans $\gg$ est égal à $\gg_0 = \aa$.

\bigskip

Posons $\alpha_i = \beta_i - \beta_{i+1}$ pour $i \in [[1,p-1]]$ et $\alpha_p = \beta_{p-1} + \beta_p$. Alors l'ensemble $(\alpha_i)_{i \in [[1,p]]}$ est une base du système de racines $\Sigma$. Les racines positives correspondantes sont
$$\Sigma^+ = \{\beta_i-\beta_j : i,j \in [[1,p]], i < j\} \cup \{\beta_i+\beta_j : i,j \in [[1,p]], i < j\}. $$
La chambre de Weyl positive (vectorielle) associée est
$$\aa^+ = \{H \in \aa : H_{1,1} > \ldots > H_{p,p} > 0 \}.$$

\bigskip

L'algèbre de Lie nilpotente $\nn = \oplus_{\alpha \in \Sigma^+} \gg_\alpha$ est égale à la sous-algèbre de Lie de $\gg$ constituée des matrices strictement triangulaires supérieures. Le sous-groupe de Lie connexe $A$ d'algèbre de Lie $\aa$ est
$$A = \{\Diag(a_1, \ldots ,a_p , a_p^{-1}, \ldots ,a_1^{-1}) : \forall i \in [[1,p]], a_i \in \;]0,+\infty[\; \}.$$
Le sous-groupe de Lie connexe $N$ d'algèbre de Lie $\nn$ est le sous-groupe de $G$ constitué des matrices triangulaires supérieures de coefficients diagonaux égaux à $1$.

\bigskip

Le sous-groupe de Lie compact maximal de $G$ d'algèbre de Lie $\kk$ est $K = G \cap \SO(2p)$. Pour voir que le groupe $K$ est isomorphe à $\SO(p) \times \SO(p)$, il est plus clair de changer la matrice de la forme quadratique de signature $(p,p)$, en prenant
$$J' = \left( \begin{array}{cc} I_p & 0 \\ 0 & -I_p \end{array} \right)$$
où $I_p$ désigne la matrice identité de $\GL(p,\R)$. Alors l'application suivante est un isomorphisme
\beq \SO(p) \times \SO(p) & \ra & K \\
(Y,Z) & \mapsto & \left( \begin{array}{cc} Y & 0 \\ 0 & Z \end{array} \right). \eeq

Le centralisateur $M=Z_K(\aa)$ de $\aa$ dans $K$ est égal à
\begin{eqnarray*} M=\{\Diag(\varepsilon_1, \ldots ,\varepsilon_p , \varepsilon_p, \ldots ,\varepsilon_1) : \forall i \in [[1,p]], \varepsilon_i = \pm 1 \mbox{ et } \Pi_{i \in [[1,p]]} \; \varepsilon_i = 1\}. \end{eqnarray*}
La condition imposée assure que les éléments considérés sont bien dans la composante neutre $G=\SO_0(p,p)$ de $\SO(p,p)$.

On remarque que, dans ce cas particulier également, les espaces de racines $\gg_\alpha$ sont tous de dimension $1$, et que le groupe $M$ est fini, mais ce n'est pas toujours le cas.

\bigskip

Le groupe $G$ est sans facteur compact, donc l'espace symétrique $G/K$ est de type non compact, de rang la dimension de $\aa$, c'est-à-dire $p$.

\section{La topologie de Chabauty sur l'espace des sous-groupes fermés}

\label{sec:chabauty}

Pour cette partie, on pourra se référer à l'excellente introduction~\cite{harpe_chabauty}.

\bigskip

Soit $X$ un espace topologique localement compact, et $\mathcal{F}(X)$ l'ensemble des fermés de $X$. On munit $\mathcal{F}(X)$ de la \df{topologie de Chabauty} (voir~\cite{chabauty}) : les ouverts sont les réunions quelconques d'intersections finies de parties de la forme :
\begin{eqnarray*}
O_K &=& \{H \in \Ch(G) : H \cap K = \emptyset \} \\
O'_U &=& \{H \in \Ch(G) : H \cap U \neq \emptyset \}
\end{eqnarray*}
où $K$ est un compact de $X$ et $U$ un ouvert de $X$.

\bigskip

Le résultat suivant est classique :

\bpro L'espace topologique $\mathcal{F}(X)$ est compact. (Voir~\cite{chabauty}, \cite[Chap.~VIII, $\S 5$]{bourbaki}, \cite[Proposition~1.7, p.~58]{cdp}.) \hfill \qed \epro

Soit $G$ un groupe topologique localement compact. On note $\Ch(G) \subset \mathcal{F}(G)$ l'ensemble de ses sous-groupes fermés, muni de la topologie induite.

\bpro \label{pro:chg_compact} Le sous-espace $\Ch(G)$ est un fermé de $\mathcal{F}(G)$, donc est compact. (Voir~\cite{chabauty}, \cite[Proposition~I.3.1.2, p.~59]{CEG}, \cite[Chap.~VIII, $\S 5$]{bourbaki}, \cite[Proposition~1.7, p.~58]{cdp}.) \hfill \qed \epro

Si $G$ est de plus muni d'une distance $d$ induisant sa topologie, on peut alors décrire la convergence d'une suite de sous-groupes fermés.

\bpro Une suite de sous-groupes fermés $(H_n)_{n \in \N}$ converge vers un sous-groupe fermé $H$ dans $\Ch(G)$ si et seulement si $H$ est l'ensemble des valeurs d'adhérence des suites de $(H_n)_{n \in \N}$, c'est-à-dire :
\begin{enumerate}
\item Pour tout $x \in H$, il existe une suite $(x_n)_{n \in \N}$ convergeant vers $x$ telle que, pour tout $n$, nous ayons $x_n \in H_n$.
\item Pour toute suite strictement croissante d'entiers $(n_k)_{k \in \N}$, pour toute suite $(x_{n_k})_{k \in \N}$ convergeant vers $x$ telle que $x_{n_k} \in H_{n_k}$ pour tout $k$, nous ayons $x \in H$.
\end{enumerate}
(Voir par exemple~\cite[Lemma~I.3.1.3, p.~60]{CEG}, \cite[Proposition~1.8, p.~60]{cdp}.) \hfill \qed \epro

\bpro \label{pro:seq} Si de plus la distance $d$ sur $G$ est propre (i.e. les boules fermées sont compactes), alors l'espace $\Ch(G)$ est métrisable, pour la distance de Hausdorff pointée. (Voir~\cite[Definition~5.43, p.~76]{bridson_haefliger}.) \hfill \qed \epro

\bpro \label{pro:ouvert} Soit $f : G \ra H$ un morphisme de groupes topologiques localement compacts, qui est une application ouverte. Alors l'application $\Ch^*(f) : \Ch(H) \ra \Ch(G)$ définie par $A \mapsto f^{-1}(A) $ est continue. \epro

\bp Soit $K$ un compact de $G$, alors $f(K)$ est un compact de $H$ et $\Ch^*(f)^{-1}(O_K) = O_{f(K)}$ est un ouvert de $\Ch(H)$. Soit $U$ un ouvert de $G$, alors par hypothèse $f(U)$ est un ouvert de $H$ et $\Ch^*(f)^{-1}(O'_U) = O'_{f(U)}$ est un ouvert de $\Ch(H)$. Ainsi l'application $\Ch^*(f)$ est continue. \ep

\bpro \label{pro:surjection} Soit $f : G \ra H$ un morphisme de groupes topologiques localement compacts, qui est une surjection ouverte. Alors l'application $\Ch^*(f) : \Ch(H) \ra \Ch(G)$ définie par $A \mapsto f^{-1}(A) $ est un homéomorphisme sur son image. \epro

\bp Puisque le morphisme $f$ est surjectif, pour tout sous-groupe fermé $A$ de $H$, nous avons $f(\Ch^*(f)(H))=H$, donc l'application $\Ch^*(f)$ est injective. D'après la proposition précédente, l'application $\Ch^*(f)$ est continue. Enfin, puisque l'espace $\Ch(H)$ est compact et l'espace $\Ch(G)$ séparé, l'application continue injective $\Ch^*(f)$ est un homéomorphisme sur son image. \ep

Voici quelques rares exemples de groupes pour lesquels le calcul explicite de l'espace des sous-groupes fermés est élémentaire.

\bpro Notons $X_\Z$ le sous-espace topologique compact de $\R$ défini par $X_\Z = \{0\} \cup \{\frac{1}{n} \,:\, n \in \N \bs \{0\} \}$. Alors l'application $\phi_\Z : X_\Z  \rightarrow \Ch(\Z)$ définie par $\frac{1}{n} \mapsto n\Z$ et $0 \mapsto \{0\}$ est un homéomorphisme. \hfill \qed \epro

\bpro \label{pro:R} L'application $\phi_\R : X_\R=[0,\infty] \rightarrow \Ch(\R)$ définie par $\alpha \mapsto \frac{1}{\alpha}\Z$ si $\alpha \in \,]0,\infty[$, $0 \mapsto \{0\}$ et $\infty \mapsto \R$ est un homéomorphisme. \hfill \qed \epro

Voici quelques exemples de groupes pour lesquels l'espace des sous-groupes fermés a été décrit : le groupe $\C$ pour lequel l'espace $\Ch(\C)$ est homéomorphe à la sphère de dimension $4$ (voir~\cite{pourezza} ou \cite{harpe_chabauty}), le groupe affine de la droite réelle (voir~\cite[Proposition~1.1, p.~2]{heisenberg}), le groupe de Heisenberg de dimension $3$ (voir~\cite[Theorem~1.3, p.~4]{heisenberg}) et le groupe $\R \times \Z$ (voir~\cite{RxZ}).

\section{La compactification de Chabauty}

Soit $X$ un espace topologique localement compact. Une \df{compactification} de $X$ est la donnée d'une paire $(K,i)$, où $K$ est un espace topologique compact et $i:X \rightarrow K$ est un plongement (i.e. $i$ réalise un homéomorphisme sur son image) d'image dense. On remarque que $i(X)$ est alors ouvert dans son adhérence. On identifie fréquemment $X$ et $i(X)$ par l'application $i$. On appelle \df{bord} de $X$ l'espace $\partial X = K \bs i(X)$. Si $H$ est un groupe agissant continûment sur $X$, on dit que c'est une \df{$H$-compactification}, ou compactification \df{$H$-équivariante}, si l'action de $H$ sur $i(X)$, conjuguée par $i$ de l'action de $H$ sur $X$, s'étend continûment à $K$. Cette extension est alors unique.

\bigskip
\textbf{Un exemple de compactification : la compactification géodésique}
\bigskip

Cette compactification peut se définir dans le cadre des espaces $\CAT(0)$ : une référence assez complète à ce sujet notamment est \cite{bridson_haefliger}.

\bigskip

Un espace métrique $X$ est dit \df{géodésique} si deux points quelconques $x$, $y$ de $X$ sont les extrémités d'un segment géodésique $[xy]$ de $X$ (non nécessairement unique, mais la notation $[xy]$ n'entraînera pas de confusion). Un espace métrique géodésique $X$ est appelé $\CAT(0)$ si, pour tout triangle géodésique $pqr$ dans $X$, et pour tous points $x \in [pq]$ et $y \in [pr]$, si l'on désigne par $\ov{p}\ov{q}\ov{r}$ le triangle euclidien de comparaison, et si $\ov{x} \in [\ov{p}\ov{q}]$ et $\ov{y} \in [\ov{p}\ov{r}]$ sont les points de comparaison dans ce triangle, alors $d(x,y) \leq d(\ov{x},\ov{y})$.


Sous cette hypothèse de courbure négative, on a alors unicité des géodésiques.

\bpro \label{pro:uniq_geod} Un espace métrique $\CAT(0)$ $X$ est uniquement géodésique, c'est-à-dire que par deux points distincts de $X$ passe une unique géodésique. (Voir~\cite[Proposition~1.4, p.~160]{bridson_haefliger}.) \hfill \qed  \epro

Voici le résultat qui permet d'affirmer que les espaces symétriques de type non compact sont des espaces métriques $\CAT(0)$.

\bthm[Cartan] Si $X$ est une variété riemannienne de courbure sectionnelle négative ou nulle, alors $X$ est un espace $\CAT(0)$. (Voir~\cite[Theorem~1A.6, p.~173]{bridson_haefliger}.) \hfill \qed  \ethm

Soit $X$ un espace métrique $\CAT(0)$ complet (en tant qu'espace métrique, c'est-à-dire tel que toute suite de Cauchy converge). Deux rayons géodésiques $c,c' : [0,\infty[\; \ra X$ sont dits \df{asymptotes} s'il existe une constante $K$ telle que, pour tout $t \in [0,\infty[$, nous ayons $d(c(t),c'(t)) \leq K$. Ceci définit une relation d'équivalence sur l'ensemble des rayons géodésiques de $X$ ; appelons \df{bord} (visuel) de $X$ l'ensemble de ses classes d'équivalence, et notons-le $\partial_\infty X$. Notons de plus $\ov{X}^g = X \cup \partial_\infty X$.

Il existe une topologie sur l'espace $\ov{X}^g$ (voir~\cite{bridson_haefliger}) qui en fait une $G$-compactification de l'espace $X$, où $G$ désigne le groupe des isométries de $X$. Elle est appelée \df{compactification géodésique}, \df{compactification conique} ou \df{compactification par le bord visuel}.

\bexe Dans le cas du plan hyperbolique réel $\H^2_\R = \SL(2,\R) / \SO(2,\R)$, dans le modèle du disque ouvert, la compactification géodésique est $\SL(2,\R)$-isomorphe à la compactification usuelle où l'on ajoute le cercle à l'infini pour obtenir le disque fermé. \eexe

Reprenons les notations de la partie~\ref{sec:espaces}, où $K$ désigne le stabilisateur dans $G$ d'un point base fixé $x_0$ de l'espace symétrique de type non compact $X$. Cet espace métrique $X$ est un espace $CAT(0)$, et on sait décrire la convergence d'une suite de points de $X$ dans la compactification géodésique. 

\bpro \label{pro:parabolique} Soit $(H_n)_{n \in \N}$ une suite de $\aa^+$ tendant vers l'infini, et soit $I = \{\alpha \in \Delta \,:\, \alpha(H_n) \ra +\infty \}$. Alors la suite $(\exp(H_n) \cdot x_0)_{n \in \N}$ de $X$ converge dans $\ov{X}^g$ vers un point de $\partial_\infty X$, dont le stabilisateur dans $G$ est le sous-groupe parabolique $P^I = K^IMAN$. (Voir~\cite[Proposition~3.9, p.~27]{guivarch}.)  \hfill \qed \epro

\bigskip
\textbf{La définition de la compactification}
\bigskip

Soit $X$ un espace symétrique de type non compact, et soit $G$ un groupe de Lie connexe muni d'une action continue isométrique sur $X$. On suppose que $G$ se surjecte sur la composante neutre $\Isom_0(X)$ du groupe des isométries de $X$, avec noyau fini. Alors le groupe de Lie $G$ est semi-simple, de centre fini et sans facteur compact. Nous verrons en fin de paragraphe que la construction qui suit est indépendante du choix d'un tel groupe $G$, ce qui permet de travailler avec des groupes de Lie ayant un centre fini non trivial, comme par exemple $\SL(n,\R)$ pour $n \geq 2$ pair.

Nous allons définir une compactification de $X$ en plongeant $X$ dans l'espace $\Ch(G)$ des sous-groupes fermés de $G$, via l'application $\phi$ qui à un point $x$ de l'espace symétrique associe son stabilisateur $G_x$ dans $G$ (lequel est fermé car l'action de $G$ sur $X$ est continue) :
\begin{eqnarray*}
\phi : X & \rightarrow & \Ch(G) \\
x & \mapsto & G_x=\{g \in G \,:\, g \cdot x=x\}.
\end{eqnarray*}

Remarquons que l'image de $\phi$ est exactement l'ensemble des sous-groupes compacts maximaux de $G$, d'après le fait que tout sous-groupe compact de $G$ fixe un point de l'espace métrique $\CAT(0)$ $X$ (voir par exemple~\cite[Proposition~2.7, p.~179]{bridson_haefliger}).

\bpro \label{pro:plongement} L'application $\phi$ est un plongement. \epro

\bp (Voir~\cite[Proposition~9.3, p.~133]{guivarch}.)

Choisissons un point base $x_0$ de l'espace symétrique $X$, et notons $K=\phi(x_0)=G_{x_0}$ le stabilisateur de $x_0$ dans $G$ : c'est un sous-groupe compact maximal. Choisissons une décomposition de Cartan $\gg = \kk \oplus \pp$, où $\kk$ est l'algèbre de Lie de $K$, une sous-algèbre de Lie abélienne maximale diagonalisable $\aa \subset \pp$ ainsi qu'une chambre de Weyl positive $\aa^+$.

\bigskip

Montrons que $\phi$ est injective : soit $g \in G$ tel que $gKg^{-1}=K$. Écrivons $g=k_1ak_2$ dans la décomposition de Cartan $K\ov{A^+}K$, alors $aKa^{-1} = K$, donc si $k \in K$ nous avons $kak^{-1} =k(aka^{-1})a \in Ka$. Or $kak^{-1} \in P$ et $a \in P$, donc par unicité de la décomposition polaire $G=KP$, nous avons $kak^{-1}=a$, et ce pour tout $k \in K$. Montrons que $a=e$.

L'exponentielle réalise un homéomorphisme de $\aa$ sur $A$ : soit donc $H \in \aa$ tel que $\exp H=a$. Pour tout $k \in K$ nous avons $aka^{-1}=k$, donc pour tout $U \in \kk$ nous avons $\Ad a(U)=e^{\ad H} U=U$. Or pour toute racine $\alpha \in \Sigma$ il existe $U_\alpha \in \gg_\alpha$ non nul, posons alors $U=U_\alpha + \theta(U_\alpha)$. Puisque $\theta(U)=U$, on en déduit que $U \in \kk$. Par ailleurs $e^{\ad H} U = e^{\alpha(H)} U_\alpha + e^{-\alpha(H)} \theta(U_\alpha)$. On ne peut donc avoir $e^{\ad H} U=U$ que si $\alpha(H)=0$, et ce pour tout $\alpha \in \Sigma$ : ainsi $H=0$. Finalement $a=e$ donc $g \in K$.

Ceci montre que $\phi$ est injective : supposons $\phi(x)=\phi(x')$. L'action de $G$ sur $X$ est transitive, soient donc $g,g' \in G$ tels que $x=g \cdot x_0$ et $x'=g' \cdot x_0$. Alors $\phi(x)=gKg^{-1}=g'Kg'^{-1}=\phi(x')$, donc $(g'^{-1}g)K(g'^{-1}g)^{-1} = K$, ainsi d'après ce qui précède $g'^{-1}g \in K$, ce qui implique $x=x'$.

\bigskip

Montrons que $\phi$ est continue. Comme toute variété topologique, le groupe de Lie $G$ est métrisable : on peut utiliser le critère séquentiel (voir la proposition~\ref{pro:seq}) pour montrer la continuité de $\phi$. Soit $g \cdot x_0 \in X$, et soit $(g_n \cdot x_0)_{n \in \N}$ une suite de $X$ convergeant vers $g \cdot x_0$. D'après la décomposition d'Iwasawa, on peut supposer que $g$ et les $g_n$ appartiennent à $AN$. Puisque $(g_n \cdot x_0)_{n \in \N}$ converge vers $g \cdot x_0$, il existe une suite $(k_n)_{n \in \N}$ dans $K$ telle que $(g_nk_n)_{n \in \N}$ converge vers $g$. D'après la continuité de la décomposition d'Iwasawa, on déduit que $(k_n)_{n \in \N}$ converge vers $e$ et que $(g_n)_{n \in \N}$ converge vers $g$. Ainsi la suite $(g_n^{-1})_{n \in \N}$ converge vers $g^{-1}$. On conclut que $(\phi(g_n \cdot x_0))_{n \in \N} = (g_nKg_n^{-1})_{n \in \N}$ converge vers $gKg^{-1}=\phi(g \cdot x_0)$ : $\phi$ est continue.

\bigskip

Montrons que l'application $\phi$ est propre : soit $(g_nKg_n^{-1})_{n \in \N}$ une suite de $\Ch(G)$ convergeant vers $gKg^{-1}$, montrons que la suite $(g_n \cdot x_0)_{n \in \N}$ de $X$ converge à extraction près.

Remarquons que, si $d$ est une distance $G$-invariante sur $X$, alors pour $h \in G$ et $k \in K$
$$d(hkh^{-1} \cdot x_0,h \cdot x_0) = d(kh^{-1} \cdot x_0,x_0) = d(h^{-1} \cdot x_0,x_0) = d(h \cdot x_0, x_0).$$
Ainsi l'orbite $hKh^{-1} \cdot x_0$ de $x_0$ sous $hKh^{-1}$ est incluse dans la sphère de centre $h \cdot x_0$, de rayon $d(h \cdot x_0,x_0)$. Donc le diamètre de $hKh^{-1} \cdot x_0$ est inférieur ou égal à $2d(h \cdot x_0,x_0)$. Soit $g \in G$ la symétrie géodésique par rapport au point $h \cdot x_0$, alors $g \in hKh^{-1}$. Or la distance de $x_0$ à $g \cdot x_0$ est égale à $2d(h \cdot x_0,x_0)$, et ces deux points sont dans l'orbite $hKh^{-1} \cdot x_0$. Par conséquent, pour tout $h \in G$, le diamètre de l'orbite $hKh^{-1} \cdot x_0$ de $x_0$ sous $hKh^{-1}$ est égal à $2d(h \cdot x_0,x_0)$.

Puisque la suite $(g_nKg_n^{-1})_{n \in \N}$ converge vers $gKg^{-1}$ dans $\Ch(G)$ et que l'évaluation en $x_0$ est continue, la suite $(g_nKg_n^{-1} \cdot x_0)_{n \in \N}$ converge vers $gKg^{-1} \cdot x_0$ dans l'espace $\mathcal{F}(M)$ des fermés de $X$, muni de la topologie de Chabauty. Cette suite appartient à partir d'un certain rang à l'ouvert $O_F$, où $F$ est le compact de $X$ défini par $F=\ov{B}(h \cdot x_0,3d(h \cdot x_0,x_0)) \bs B(h \cdot x_0,2d(h \cdot x_0,x_0))$. De plus, le complémentaire de $F$ dans $X$ a deux composantes connexes, $M \bs \ov{B}(h \cdot x_0,3d(h \cdot x_0,x_0))$ et $B(h \cdot x_0,2d(h \cdot x_0,x_0))$. Or à partir d'un certain rang, l'ensemble $g_nKg_n^{-1} \cdot x_0$ intersecte $B(h \cdot x_0,2d(h \cdot x_0,x_0))$, donc par connexité de $g_nKg_n^{-1} \cdot x_0$ (car $K$ est connexe), l'orbite $g_nKg_n^{-1} \cdot x_0$ est incluse dans $B(h \cdot x_0,2d(h \cdot x_0,x_0))$ à partir d'un certain rang. Ainsi le diamètre de $g_nKg_n^{-1} \cdot x_0$, qui est égal à $2d(g_n \cdot x_0,x_0)$, est borné. Donc la suite $(g_nK)_{n \in \N}$ est bornée : elle converge, à extraction près.

\bigskip

L'application $\phi$, continue, injective et propre, est donc un homéomorphisme sur son image.
\ep

Notons $\ov{X}^\Ch$ l'adhérence de $\phi(X)$ dans $\Ch(G)$ : l'espace $\ov{X}^\Ch$ est compact. Le couple $(\ov{X}^\Ch,\phi)$ est appelé la \df{compactification de Chabauty} de l'espace symétrique $X$. Le groupe $G$ agit sur $\Ch(G)$ par conjugaison, et l'application $\phi$ est $G$-équivariante : si $g \in G$ et $x \in X$, le stabilisateur de $g \cdot x$ dans $X$ est $G_{g \cdot x} = gG_xg^{-1} = g \cdot G_x$. Ainsi, la compactification de Chabauty $\ov{X}^\Ch$ est une $G$-compactification de $X$.

\bigskip

Notons $G^0=\Isom_0(X)$ la composante neutre du groupe des isométries de $X$. Nous allons montrer que la compactification de Chabauty de $X$ dans $\Ch(G)$ ne dépend pas, à isomorphisme près, du choix du groupe de Lie connexe $G$ agissant continûment isométriquement sur $X$, et se surjectant avec noyau fini sur $G^0$.

\bigskip

Notons $\pi : G \ra G^0$ la projection définie par l'action de $G$ sur $X$. Par hypothèse, $\pi$ est un morphisme de groupes de Lie surjectif de noyau fini. En particulier, $\pi$ est un difféomorphisme local. Soit $K^0$ un sous-groupe compact maximal de $G^0$, et $K$ le sous-groupe compact maximal $K=\pi^{-1}(K^0)$ de $G$. Choisissons une métrique riemannienne $G^0$-invariante à gauche, $K^0$-invariante à droite sur $G^0$. Alors la métrique riemannienne tirée en arrière sur $G$ par $\pi$ est $G$-invariante à gauche et $K$-invariante à droite.

Notons $\phi : X \ra \Ch(G)$ (resp. $\phi^0 : X \ra \Ch(G^0)$) l'application qui à un point $x \in X$ associe son stabilisateur dans $G$ (resp. dans $G^0$) : ce sont deux plongements de l'espace symétrique $X$. Notons $\ov{X}^\Ch$ (resp. $\ov{X}^{\Ch,0}$) la compactification associée.

D'après la proposition~\ref{pro:surjection}, la projection $\pi : G \ra G^0$ induit un plongement $\Ch^*(\pi) : \Ch(G^0) \ra \Ch(G)$.

\bpro L'application $\Ch^*(\pi) : \Ch(G^0) \ra \Ch(G)$ est un homéomorphisme de $\ov{X}^{\Ch,0}$ sur $\ov{X}^\Ch$ qui \df{entrelace} $\pi$, i.e. pour tous $g \in G$ et $x \in \ov{X}^{\Ch,0}$ nous avons $\Ch^*(\pi)(\pi(g) \cdot x) = g \cdot \Ch^*(\pi)(x)$. \epro

\bp
On constate que $\phi = \Ch^*(\pi) \circ \phi^0$, donc ceci montre que $\Ch^*(\pi)(\ov{X}^{\Ch,0}) = \ov{X}^\Ch$. Or l'application $\Ch^*(\pi)$ est un plongement, donc c'est un homéomorphisme de $\ov{X}^{\Ch,0}$ sur $\ov{X}^\Ch$.

De plus, pour tous $g \in G$ et $x \in \ov{X}^{\Ch,0} \subset \Ch(G^0)$ nous avons $\pi(g) \cdot x = \pi(g)x\pi(g)^{-1} = \pi(gxg^{-1}) = \pi(g \cdot x)$. Ainsi $g \cdot x = \Ch^*(\pi) (\pi(g) \cdot x)$.
\ep

Ceci montre par exemple que, pour étudier la compactification de Chabauty de l'espace symétrique associé à $\SL(n,\R)$, on peut travailler dans l'espace $\Ch(\SL(n,\R))$, plus maniable que $\Ch(\PSL(n,\R))$.

\bigskip
\textbf{La distalité}
\bigskip

Si $X$ est un espace métrique et $G$ est un groupe d'homéomorphismes de $X$, on dit que l'action de $G$ sur $X$ est \df{distale} si pour tous $x,y \in X$ distincts il existe $\varepsilon >0$ tel que pour tout $g \in G$ nous ayons $d(g \cdot x,g \cdot y) \geq \varepsilon$. Si $g$ est un homéomorphisme de $X$, on dit que son action sur $X$ est distale si l'action du groupe $\{g^n \,:\, n \in \Z\}$ engendré par $g$ est distale.

\bigskip

Si $V$ est un espace vectoriel réel de dimension finie, on dit qu'un élément $h \in \GL(V)$ est
\df{distal} si son action linéaire sur $V$ est distale, ce qui équivaut à demander que le spectre de $h$ soit inclus dans le cercle unité $\SS^1$. D'après le théorème~1 de \cite{conze_guivarch}, l'action linéaire d'un groupe $H$ sur $V$ est distale si et seulement si tout élément de $H$ est distal.

\bigskip

Un sous-groupe $H$ d'un groupe de Lie $G$ est appelé un \df{sous-groupe distal} de $G$ si son action adjointe sur l'algèbre de Lie de $G$ est distale, ce qui équivaut à demander que pour tout $h \in H$ le spectre de $\Ad h$ agissant sur l'algèbre de Lie de $G$ soit inclus dans $\SS^1$.

\bigskip

Cette notion, introduite par Furstenberg, provient des systèmes dynamiques (voir par exem\-ple~\cite[Proposition~9.5, p.~135]{guivarch}). L'intérêt d'introduire ici la distalité réside dans la proposition suivante.

\bpro \label{pro:distal_limite} Soit $G$ un sous-groupe de $\GL(V)$, où $V$ est un espace vectoriel réel de dimension finie. Le sous-espace $\Ch_{\mbox{\scriptsize distal}} (G)$ de $\Ch(G)$, constitué des sous-groupes fermés de $G$ dont l'action sur $V$ est distale, est fermé. \epro

\bp Soit $H$ un sous-groupe fermé de $G$ adhérent à $\Ch_{\mbox{\scriptsize distal}}(G)$, et soit $h \in H$. Soit $(h_n)_{n \in \N}$ une suite d'éléments de $G$ ayant une action distale sur $V$, convergeant vers $h$. On peut voir le spectre comme application de $\GL(V)$ à valeurs dans l'ensemble des $n$-uplets de nombres complexes, à permutation près. On munit cet espace, naturellement identifié à $\C^n / \mathfrak{S}_n$, de la topologie quotient de $\C^n$ par l'action du groupe symétrique $\mathfrak{S}_n$ par permutation. Alors le spectre est continu. Puisque le cercle unité $\SS^1$ est fermé, le spectre de $h$ est lui aussi inclus dans $\SS^1$. Ainsi, $H$ est un sous-groupe de $G$ ayant une action distale sur $V$. \ep

\bpro \label{pro:compact_distal} Tout sous-groupe de Lie compact d'un groupe de Lie est distal. \epro

\bp Soit $H$ un sous-groupe de Lie compact d'un groupe de Lie $G$. Soit $g \in H$, alors la suite $(g^n)_{n \in \N}$ admet une valeur d'adhérence $h$ ; et si $\lambda$ est une valeur propre de $\Ad g$, alors $\lambda^n$ converge à extraction près vers une valeur propre de $\Ad h$. Or $\Ad h$ est inversible, ce qui implique que $|\lambda| = 1$. Donc, pour tout $g \in H$, le spectre de $\Ad g$ est inclus dans le cercle unité $\SS^1$ : ainsi, $H$ est un sous-groupe distal de $G$. \ep

\bpro \label{pro:distal} La compactification de Chabauty $\ov{X}^\Ch \subset \Ch(G)$ de $X$ est constituée de sous-groupes distaux de $G$. \epro

\bp Par définition de la compactification de Chabauty, l'espace $\ov{X}^\Ch$ est l'adhérence des sous-groupes compacts maximaux de $G$, donc d'après les propositions~\ref{pro:distal_limite} et \ref{pro:compact_distal}, l'espace $\ov{X}^\Ch$ est constitué de sous-groupes distaux de $G$. \ep

Concernant l'étude des actions distales de groupes, on pourra se référer à~\cite{abels_affine}, \cite{abels_lie} et \cite{conze_guivarch}.

\bigskip
\textbf{Détermination des groupes limites}
\bigskip

Fixons $I$ une partie propre de la base $\Delta$ du système de racines $\Sigma$. L'ensemble $\ov{\aa_I^+} = \aa_I \cap \ov{\aa^+}$ est muni d'un structure de cône simplicial fermé :
$$\ov{\aa_I^+} = \left( \bigcap_{\alpha \in I} \Ker \alpha \right) \cap \left( \bigcap_{\alpha \in \Delta \bs I} \{H \in \aa : \alpha(H) \geq 0\} \right).$$

Les facettes du cône $\ov{\aa_I^+}$ sont constituées des
$$\left( \bigcap_{\alpha \in J} \Ker \alpha \right) \cap \left( \bigcap_{\alpha \in \Delta \bs J} \{H \in \aa : \alpha(H) \geq 0\} \right),$$
où $J$ est une partie de $\Delta$ contenant strictement $I$.

Notons $\aa_I^+$ l'intérieur de $\ov{\aa_I^+}$ dans l'espace vectoriel $\aa_I$ : c'est le cône simplicial ouvert
$$\aa_I^+ = \left( \bigcap_{\alpha \in I} \Ker \alpha \right) \cap \left( \bigcap_{\alpha \in \Delta \bs I} \{H \in \aa : \alpha(H) > 0\} \right).$$

\bigskip

On dit qu'une suite $(H_n)_{n \in \N}$ d'éléments de $\aa_I^+$ \df{tend vers $+\infty$ dans $\aa_I^+$} si, pour toute racine $\alpha \in \Delta \bs I$, nous avons $\alpha(H_n) \ral{n \ra +\infty} +\infty$, c'est-à-dire que la suite $(H_n)_{n \in \N}$ s'éloigne de chacun des murs du cône $\aa_I^+$. L'exponentielle réalisant un difféomorphisme de $\aa_I^+$ sur $A_I^+=A_I \cap A^+$, on dit qu'une suite $(a_n)_{n \in \N}$ d'éléments de $A_I^+$ \df{tend vers $+\infty$ dans $A_I^+$} si la suite $(\log a_n)_{n \in \N}$ tend vers $+\infty$ dans $\aa_I^+$.

\bigskip

Notons de plus $D^I = K^IMN_I$.

\bexe Pour $G=\SL(n,\R)$, $K=\SO(n,\R)$ et $I = \{\alpha_{i,i+1},\alpha_{j,j+1}\}$, où $1 \leq i < j \leq n-1$ :
$$D^I = \left\{ \left( \begin{array}{ccc} U_i & * & * \\0 & U_{j-i} & * \\ 0 & 0 & U_{n-j} \end{array} \right) \in G : U_k \in \O(k,\R) \right\}.$$  
\eexe

\bpro \label{pro:kimni} Le sous-groupe $M$ normalise $K^I$, donc $K^IM$ est un sous-groupe de Lie de $G$. Le sous-groupe $K^IM$ normalise $N_I$, donc $D^I = K^IMN_I$ est un sous-groupe de Lie de $G$, d'algèbre de Lie notée $\ddd^I = (\kk^I+\mm) \oplus \nn_I$. \epro

\bp Montrons que le sous-groupe $M$ normalise l'algèbre de Lie $\zz_\gg(\aa_I)$ : soient $m \in M$, $X \in \zz_\gg(\aa_I)$ et $H \in \aa_I$. Puisque $M$ centralise $\aa$, nous avons
$$\ad (\Ad m (X)) H = [\Ad m (X), H] = \Ad m ([X,\Ad m^{-1} (H)]) = \Ad m ([X,H]) = 0.$$
Ainsi $\Ad m (X) \in \zz_\gg(\aa_I)$ donc $M$ normalise $\zz_\gg(\aa_I)$. On en déduit que $M$ normalise également son algèbre dérivée $D(\zz_\gg(\aa_I)) = \gg^I$, puis que $M$ normalise le sous-groupe de Lie connexe $G^I$ d'algèbre de Lie $\gg^I$. Or $M$ est un sous-groupe de $K$, donc $M$ normalise $G^I \cap K = K^I$. Ainsi, $K^IM$ est un sous-groupe compact de $G$.

\bigskip

Montrons que le groupe $M$ normalise l'algèbre de Lie $\nn_I$ : soient $m \in M$, $\alpha \in \Sigma$, $Y \in \gg_\alpha$ et $H \in \aa$. Alors
$$\ad (\Ad m (Y)) H = \Ad m ([Y,\Ad m^{-1} (H)]) = \Ad m ([Y,H]) = \alpha(H) \Ad m (Y).$$
Donc le groupe $M$ normalise tous les espaces de racines $\gg_\alpha$, donc en particulier $M$ normalise l'algèbre de Lie $\nn_I$.

Montrons que le groupe $K^I$ normalise l'algèbre de Lie $\nn_I$ : soit $Y \in \kk^I$, écrivons la décomposition de $Y$ en espaces de racines $Y = Y_0 + \Sigma_{\alpha \in \Sigma^I} Y_\alpha$, où $Y_0 \in \gg_0$ et $Y_\alpha \in \gg_\alpha$ pour toute racine $\alpha \in \Sigma^I$. Soit $Z \in \nn_I$, écrivons sa décomposition en espaces de racines $Z = \Sigma_{\beta \in \Sigma_I^+} Z_\beta$, où $Z_\beta \in \gg_\beta$ pour toute racine $\beta \in \Sigma_I^+$. Soient $\alpha \in \Sigma^I$ et $\beta \in \Sigma_I^+$. Si $\alpha+\beta$ est une racine, celle-ci doit appartenir à $\Sigma_I^+$, car dans l'écriture de $\alpha+\beta$ sur la base $\Delta$, l'une des coordonnées des racines de $\Delta \bs I$ est strictement positive. Donc $[Y,Z] \in \nn_I$. Ainsi l'algèbre de Lie $\kk^I$ normalise $\nn_I$, puis le groupe connexe $K^I$ normalise $\nn_I$.

En conclusion, le groupe $K^IM$ normalise l'algèbre de Lie $\nn_I$ et également le groupe de Lie $N_I$. Par conséquent $D^I=K^IMN_I$ est un sous-groupe de Lie de $G$ (le groupe $D^I$ est fermé car $K^IM$ est compact et $N_I$ est fermé). \ep

Nous allons maintenant étudier les \df{groupes limites}, c'est-à-dire les éléments de $\ov{X}^\Ch \bs \phi(X)$.

\bpro \label{pro:di} Soit $(a_n)_{n \in \N}$ une suite de $A_I^+$ tendant vers $+\infty$ dans $A_I^+$. Alors la suite de sous-groupes fermés $(a_nKa_n^{-1})_{n \in \N}$ converge vers $D^I$ dans $\Ch(G)$. \epro

\bexe Pour $G=\SL(3,\R)$, $K=\SO(3,\R)$ et $I=\{\alpha_{1,2}\}$, soit
$$ a_n = \left( \begin{array}{ccc} \lambda_n & 0 & 0 \\0 & \lambda_n & 0 \\ 0 & 0 & \frac{1}{\lambda_n^2} \end{array} \right) \in A_I^+, $$
où $\lambda_n \ral{n \ra +\infty} +\infty$.
Alors on peut vérifier dans ce cas particulier que la suite de sous-groupes fermés $a_n \SO(3,\R) a_n^{-1}$ converge vers le sous-groupe fermé
$$D^I = \left( \begin{array}{cc} \O(2) & * \\ 0 & \pm 1 \end{array} \right) \subset \SO(3,\R).$$
\eexe 

Fixons $(a_n)_{n \in \N}$ comme dans cet énoncé. Puisque $\Ch(G)$ est compact, pour montrer cette proposition, il suffit de montrer que toute valeur d'adhérence de la suite $(a_nKa_n^{-1})_{n \in \N}$ est égale à $D^I$ : soit $D$ une valeur d'adhérence. Supposons, quitte à extraire une sous-suite, que la suite $(a_nKa_n^{-1})_{n \in \N}$ converge vers $D$. La preuve de cette proposition, que nous donnerons plus loin, découle immédiatement des lemmes suivants.

\blem \label{lem:inclusion} Nous avons l'inclusion $D^I \subset D$. \elem

\bp Montrons que $N_I \subset D$ : soit $y \in N_I$. Notons $H_n \in \aa_I$ et $Y \in \nn_I$ tels que $\exp H_n=a_n$ et $\exp Y=y$. \'{E}crivons de plus $Y=\Sigma_{\alpha \in \Sigma_I^+}Y_\alpha$, où $Y_\alpha \in \gg_\alpha$ pour tout $\alpha \in \Sigma_I^+$. Alors
$$a_n^{-1}ya_n = \exp( \Ad a_n^{-1}(Y)) = \exp( \Sigma_{\alpha \in \Sigma_I^+} e^{-\ad H_n} Y_\alpha) = \exp( \Sigma_{\alpha \in \Sigma_I^+} e^{-\alpha(H_n)} Y_\alpha).$$
Or l'hypothèse que $a_n$ tend vers $+\infty$ dans $A_I^+$ signifie que $\alpha(H_n) \ral{n \rightarrow +\infty} +\infty$ pour tout $\alpha \in \Sigma_I^+$. Donc $a_n^{-1}ya_n \ral{n \rightarrow +\infty} e$.
Or, d'après la décomposition d'Iwasawa $G=KAN^-$, soient $k_n \in K$, $a'_n \in A$ et $y'_n \in N^-$ tels que $a_n^{-1}ya_n = k_n a'_n y'_n$. Par continuité de la décomposition d'Iwasawa, nous avons $k_n \ral{n \rightarrow +\infty} e$, $a'_n \ral{n \rightarrow +\infty} e$ et $y'_n \ral{n \rightarrow +\infty} e$. Puisque $N^-$ est égal au produit semi-direct $N^- = N_I^- \rtimes N^{I,-}$, écrivons $y'_n=y'^I_ny'_{n,I}$, où $y'^I_n \in N^{I,-}$ et $y'_{n,I} \in N^-_I$ : on a alors $y'^I_n \ral{n \rightarrow +\infty} e$ et $y'_{n,I} \ral{n \rightarrow +\infty} e$.

Ainsi, puisque $N^{I,-}$ et $A_I$ commutent, nous avons $a_ny'_na_n^{-1}=y'^I_na_ny'_{n,I}a_n^{-1}$. Soit $Y'_{n,I} \in \nn_I^-$ tel que $\exp Y'_{n,I}=y'_{n,I}$, alors on peut écrire $Y'_{n,I}=\Sigma_{\alpha \in \Sigma_I^-} \; Y'_{n,\alpha}$, où $Y'_{n,\alpha} \in \gg_\alpha$ pour tout $\alpha \in \Sigma_I^-$. Ainsi, $a_ny'_{n,I}a_n^{-1} = \exp( \Sigma_{\alpha \in \Sigma_I^-} \; e^{\alpha(H_n)} Y'_{n,\alpha})$. Comme $\alpha \in \Sigma_I^-$ et $H_n \in \aa_I$, nous avons $\alpha(H_n) \ral{n \rightarrow +\infty} -\infty$, donc $a_ny'_na_n^{-1} \ral{n \rightarrow +\infty} e$. Puisque $a'_n \ral{n \rightarrow +\infty} e$, on en déduit que $y = (a_nk_na_n^{-1})a'_n(a_ny'_na_n^{-1}) = \liml_{n \rightarrow +\infty} a_nk_na_n^{-1} \in D$.

\bigskip

Montrons que $K^IM \subset D$. Puisque $M=Z_K(\aa)$, on en déduit que $M$ commute à $A_I$. Par ailleurs, puisque $K^I \subset G^I \subset Z_G(\aa_I)$, le groupe $K^I$ commute à $A_I$. Soit $k \in K^IM \subset K$, alors $k=a_nka_n^{-1}$ pour tout entier $n \in \N$, donc $k \in D$.

Finalement $D^I = K^IMN_I \subset D$.
\ep 

\blem \label{lem:infini} Nous avons l'inclusion $D \subset P^I$. \elem

\bp Plaçons-nous dans la compactification géodésique $\ov{X}^g$ de $X$. Puisque la suite $(a_n)_{n \in \N}$ tend vers l'infini dans $A_I^+$, d'après la proposition~\ref{pro:parabolique}, la suite $(x_n = a_n \cdot x_0)_{n \in \N}$ converge vers un point $x \in \partial_\infty X$, dont le stabilisateur dans $G$ est $P^I$. Soit $d \in D$, et soit $(k_n)_{n \in \N}$ une suite de $K$ telle que la suite $(a_nk_na_n^{-1})_{n \in \N}$ converge vers $d$. Puisque l'action de $G$ sur la compactification $\ov{X}^g$ est continue, on en déduit que la suite $(a_nk_na_n^{-1} \cdot x_n)_{n \in \N} =  (x_n)_{n \in \N}$ converge vers $d \cdot x = x$. Par conséquent, l'élément $d \in D$ fixe le point à l'infini $x$, donc appartient au sous-groupe $P^I$. \ep

\blem \label{lem:inclusion_reciproque} Nous avons l'inclusion $D \subset D^I$. \elem

\bp Soit $d \in D$, montrons que $d \in D^I$. Puisque $d \in P^I = K^IMAN$ et que $N=N^IN_I$, on peut supposer que $d \in AN^I$. D'après la proposition~\ref{pro:distal}, l'élément $d$ est distal. Or cet élément s'écrit $d=au$, où $a \in A$ et $u \in N^I$. Soit $H \in \aa$ tel que $\exp H = a$ : alors, pour toute racine $\alpha \in \Sigma$, le réel $\exp(\alpha(H))$ est une valeur propre de $\Ad d$. Par distalité, on en déduit que pour tout $\alpha \in \Sigma$ nous avons $\alpha(H)=0$ : ainsi $H=0$ et $d=u \in N^I$. Puisque l'exponentielle du groupe de Lie $N$ est un difféomorphisme, soit $Y \in \nn^I$ tel que $\exp Y = d$.

\bigskip

Soit $(k_n)_{n \in \N}$ une suite de $K$ telle que
$$a_nk_na_n^{-1} \ral{n \ra +\infty} \exp Y .$$
Soit $H_n \in \aa_I$ tel que $\exp H_n=a_n$. Puisque le spectre de $\ad Y \in \gl(\gg)$ est égal à $\{0\}$, d'après~\cite[Theorem~1.7, p.~105]{helgason}, l'exponentielle est un difféomorphisme local en $Y$. Donc il existe une suite $(Y_n)_{n \in \N}$ de $\gg$ convergeant vers $Y$ telle que $\exp Y_n=a_nk_na_n^{-1}$ pour tout entier $n \in \N$. Posons $U_n = \Ad (a_n^{-1}) (Y_n) \in \gg$, de sorte que $\exp U_n = k_n \in K$. Notons $Y'_n = \Ad a_n\, \theta(U_n)$, de sorte que $\exp Y'_n = \exp Y_n$ : nous allons montrer que $Y'_n=Y_n$.

\bigskip

On remarque que pour tout $\alpha \in \Sigma^{I,+}$ nous avons $\alpha(H_n)=0$ donc
$$Y'_{n,\alpha} = \theta(U_n)_\alpha = \theta(U_{n,-\alpha}) = \theta(Y_{n,-\alpha}) \ral{n \ra +\infty} \theta(Y_{-\alpha}) = 0.$$
Et pour tout $\alpha \in \Sigma_I^-$ nous avons $\alpha(H_n) \ral{n \ra +\infty} -\infty$ et
$$ e^{-\alpha(H_n)}\theta(U_{n,-\alpha}) = \theta(Y_{n,-\alpha}) \ral{n \ra +\infty} \theta(Y_{-\alpha}) = 0, $$
$$\mbox{donc } Y'_{n,\alpha} = e^{\alpha(H_n)}\theta(U_{n,-\alpha}) \ral{n \ra +\infty} 0.$$
Enfin, la composante $Y'_{n,0}$ de $Y'_n$ sur $\gg_0$ vérifie
$$ Y'_{n,0} = \theta(U_{n,0}) = \theta(Y_{n,0}) \ral{n \ra +\infty} \theta(Y_0) = 0.$$

Ainsi la distance de $Y'_n$ à la sous-algèbre nilpotente $\nn^{I,-} \oplus \nn_I$ tend vers $0$, donc le spectre de la suite d'endomorphismes $(\ad Y'_n)_{n \in \N}$ converge vers $\{0\}$. Si on appelle $\mathcal{E} = \{Z \in \gg \,:\, |\Im \Sp(\ad Z)| < \pi\}$, cela signifie que pour $n$ assez grand nous avons $Y'_n \in \mathcal{E}$. Par ailleurs $Y \in \nn^I$ donc $Y_n \in \mathcal{E}$ pour $n$ assez grand. Or d'après~\cite[Théorème~3.8.4, p.~83]{mneimne_testard}, l'exponentielle de $G$ est un difféomorphisme de $\mathcal{E}$ sur son image, donc comme nous avons $\exp Y_n=\exp Y'_n$, on en déduit que $Y_n=Y'_n$ si $n$ est assez grand.

\bigskip 

Or les valeurs d'adhérence de la suite $(Y'_n)_{n \in \N}$ appartiennent à $\nn^{I,-} \oplus \nn_I$, donc $Y \in (\nn^{I,-} \oplus \nn_I) \cap \nn^I = \{0\}$. On a donc montré l'inclusion $D \subset D^I$. \ep

Notons, pour tous $I \subset \Delta$, $a \in A$ et $k \in K$,
$$D^I_{a,k} = kaD^Ia^{-1}k^{-1} = kaK^IMN_Ia^{-1}k^{-1}.$$

\bthm \label{thm:existence_ecriture} Soit $D \in \ov{X}^\Ch \bs X$ un groupe limite. Alors il existe $I$ une partie propre de $\Delta$, $a \in \ov{A^{I,+}}$ et $k \in K$ tels que $D=D^I_{a,k}$. \ethm

\bp Soit $(g_n)_{n \in \N}$ une suite de $G$ telle que $D = \liml_{n \ra +\infty} g_nKg_n^{-1}$. Soit $g_n=k_n a_n k'_n$ l'écriture de $g_n$ dans la décomposition $K\ov{A^+}K$ de $G$ : ainsi,  $k_n,k_n' \in K$ et $a_n \in \ov{A^+}$. Soit $I \subset \Delta$ l'ensemble des $\alpha \in \Delta$ tels que la suite $(\alpha(\log a_n))_{n \in \N}$ soit bornée. On ne peut avoir $I=\Delta$ car sinon la suite $(g_n)_{n \in \N}$ serait bornée, donc convergerait à extraction près vers un élément $g \in G$, or $D=gKg^{-1}$ n'est pas un groupe limite.

\bigskip

Quitte à extraire une sous-suite, on peut supposer que la suite $(k_n)_{n \in \N}$ converge vers $k \in K$ et que :
\beq & \forall \alpha \in I,& \mbox{la suite }(\alpha(\log a_n))_{n \in \N} \mbox{ converge} \\
\mbox{et } & \forall \alpha \in \Delta \bs I,& \mbox{la suite }(\alpha(\log a_n))_{n \in \N} \mbox{ tend vers }+\infty .
\eeq
Alors, si l'on écrit la décomposition $a_n = a_n^I a_{n,I}$, où $a_n^I \in A^I$ et $a_{n,I} \in A_I$, ceci revient à dire que $(a_n^I)_{n \in \N}$ converge vers $a \in A^I \cap \ov{A^+} = \ov{A^{I,+}}$ et que $(a_{n,I})_{n \in \N}$ tend vers $+\infty$ dans $A_I$. D'après la proposition~\ref{pro:di}, la suite $a_{n,I}Ka_{n,I}^{-1}$ converge dans $\Ch(G)$ vers $D^I$. Donc
$$D = \liml_{n \ra +\infty} g_nKg_n^{-1} = \liml_{n \ra +\infty} k_n a_n^I (a_{n,I}K(a_{n,I})^{-1}) {a_n^I}^{-1} k_n^{-1} = k a D^I a^{-1} k^{-1} = D^I_{a,k}. $$
\ep

Nous démontrons ensuite la fin du théorème annoncé dans l'introduction.

\bthm \label{thm:unicite_ecriture} Cette écriture est unique au sens suivant : soient $I_1,I_2$ deux parties propres de $\Delta$, $a_1 \in \ov{A^{I_1,+}}$, $a_2 \in \ov{A^{I_2,+}}$ et $k_1,k_2 \in K$. Alors nous avons $D^{I_1}_{a_1,k_1}=D^{I_2}_{a_2,k_2}$ si et seulement si
$$ I_1=I_2=I, a_1 = a_2 =a \mbox{ et } k_2^{-1}k_1 \in (K^I \cap aK^Ia^{-1})M. $$
\ethm

\blem \label{lem:ru_di} L'ensemble des vecteurs $X$ de $\ddd^I$ tels que $\ad_\gg X$ soit un endomorphisme nilpotent de $\gg$ est la sous-algèbre de Lie $\nn_I$. \elem

\bp Pour tout $X \in \nn_I \subset \nn$, d'après la décomposition de $\gg$ en espaces de racines, on sait que l'endomorphisme $\ad_\gg X$ est nilpotent. Réciproquement, soit $X \in \ddd^I$ tel que $\ad_\gg X$ soit nilpotent : écrivons $X= U + Y$, où $U \in \kk^I + \mm$ et $Y \in \nn_I$. Puisque $[\kk^I+\mm,\nn_I] \subset \nn_I$, pour tout entier $n \in \N$ nous avons $(\ad_\gg X)^n = (\ad_\gg U)^n + \ad_\gg Y_n$, où $Y_n \in \nn_I$. Or il existe un entier $n \in \N$ tel que $(\ad_\gg X)^n=0$, donc tel que $(\ad_\gg U)^n = -\ad_\gg Y_n$ soit nilpotent : ainsi l'endomorphisme $\ad_\gg U$ est nilpotent, et antisymétrique pour le produit scalaire $B_\theta$ puisque $U \in \kk$. Donc $\ad_\gg U=0$, puis $U \in \zz(\gg)=\{0\}$ car l'algèbre de Lie $\gg$ est semi-simple. Ainsi $X \in \nn_I$.
\ep

\blem \label{lem:norm_ni} Le normalisateur de l'algèbre de Lie $\nn_I$ dans $G$ est $P^I=K^IMAN$. \elem

\bp L'algèbre de Lie $\nn_I$ est normalisée par le sous-groupe $AN$.

Soit $k \in K$, et soit $U \in \kk$ tel que $\exp U = k$. Écrivons la décomposition $U = U_0 + \Sigma_{\alpha \in \Sigma} U_\alpha$ de $U$ en espaces de racines.

Si $k \in K^IM$, alors $U \in \gg_0 + \gg^I$. Pour toutes les racines $\alpha \in \Sigma^I$ et $\beta \in \Sigma_I^+$, si $\alpha+\beta$ est une racine, nous avons $\alpha + \beta \in \Sigma_I^+$, donc $[U_\alpha,\gg_\beta] \subset \nn_I$. Par ailleurs $[U_0,\gg_\beta] \subset \gg_\beta \subset \nn_I$. Par conséquent, nous avons $[U,\nn_I] \subset \nn_I$. Ainsi $k$ normalise l'algèbre de Lie $\nn_I$.

Si $k \not\in K^IM$, puisque $\Lie(K^IM)=(\mm + \gg^I) \cap \kk = \left(\gg_0 \oplus \bigoplus_{\alpha \in \Sigma^I} \gg_\alpha \right)$, il existe une racine $\alpha \in \Sigma_I^+$ telle que $U_\alpha \neq 0$. Alors $U_\alpha \in \nn_I$, et la composante sur $\gg_0$ de $[U,U_\alpha]$ dans la décomposition en espaces de racines est $[U_{-\alpha},U_\alpha] = [\theta (U_\alpha),U_\alpha]$ car $U \in \kk$. Soit $H \in \aa$ tel que $\alpha(H) \neq 0$. Alors
\beq B([\theta (U_\alpha),U_\alpha],H) &=& B([U_\alpha,H],\theta(U_\alpha)) \\
 &=& B(-\alpha(H)U_\alpha,\theta(U_\alpha)) \\
 &=& \alpha(H) B_\theta(U_\alpha,U_\alpha) \neq 0 \eeq
Ainsi, $[\theta (U_\alpha),U_\alpha] \neq 0$ : par conséquent, le vecteur $[U,U_\alpha]$ n'appartient pas à $\nn_I$, donc $U$ n'appartient pas au normalisateur de $\nn_I$, et par conséquent $k$ ne normalise pas l'algèbre de Lie $\nn_I$.

En résumé, le normalisateur de l'algèbre de Lie $\nn_I$ dans $G$ est $K^IMAN$. \ep

\blem \label{lem:norm_di} Le normalisateur du sous-groupe $D^I$ dans $G$ est $K^IMA_IN_I$. \elem

\bp Notons $L$ le normalisateur du sous-groupe $D^I$ dans $G$. Alors $L$ normalise éga\-lement la sous-algèbre $\nn_I$ des éléments $\ad_\gg$-nilpotents de son algèbre de Lie $\ddd^I$. Ainsi $L$ est inclus dans le normalisateur $P^I = K^IMAN$ de l'algèbre de Lie $\nn_I$. Soit $x \in A^IN^I \cap L$, alors pour tout $k \in K^I$, nous avons $xkx^{-1} \in D^I \cap G^I$ d'après la décomposition d'Iwasawa $G^I=K^IA^IN^I$. Or $D^I \cap G^I = K^I$, donc $xkx^{-1} \in K^I$. Ceci étant vrai pour tout $k \in K^I$, on en déduit que $xK^Ix^{-1} = K^I$. D'après la proposition~\ref{pro:plongement} appliquée au groupe $G^I$, on en déduit que $x \in K^I$. Par conséquent, le normalisateur du sous-groupe $D^I$ dans $G$ est $K^IMA_IN_I$. \ep

\bpthm (voir~\cite[Proposition~9.116, p.~139]{guivarch})

Si les groupes $k_1a_1 D^{I_1} {a_1}^{-1} k_1^{-1}$ et $k_2a_2 D^{I_2} {a_2}^{-1} k_2^{-1}$ sont égaux, on en déduit que les sous-algèbres constituées d'éléments $\ad_\gg$-nilpotents de leurs algèbres de Lie sont égales, donc d'après le lemme~\ref{lem:ru_di} nous avons $\Ad(k_1a_1) \nn_{I_1} = \Ad(k_2a_2) \nn_{I_2}$. Or $A$ normalise $\nn_{I_1}$ et $\nn_{I_2}$, donc si l'on pose $k = k_1^{-1}k_2$ nous avons $\nn_{I_1} = \Ad(k) \nn_{I_2}$. Par conséquent, nous avons également $N_{I_1} = kN_{I_2}k^{-1}$, donc d'après le lemme~\ref{lem:norm_ni} leurs normalisateurs $P^{I_1} = kP^{I_2}k^{-1}$ sont égaux.

Or, d'après le paragraphe~\cite[1.2.3, pp.~55-70]{warner}, la paire $(G,MAN)$ forme un système de Tits, donc d'après le théorème~\cite[Theorem~1.2.1.1, p.~46]{warner}, les sous-groupes paraboliques $(P^I)_{I \subset \Delta}$ sont deux à deux non conjugués, et chacun est égal à son propre normalisateur. Ainsi, $I_1=I_2=I$ et $k \in P^I$. Or $k \in K$, donc $k \in P^I \cap K = K^IM$.

Enfin $D^I = (a_1^{-1}ka_2) D^I (a_1^{-1}ka_2)^{-1}$, donc d'après le lemme~\ref{lem:norm_di}, nous avons $a_1^{-1}ka_2 \in K^IMA_IN_I$. Or $a_1^{-1}ka_2 \in A^IK^IMA^I \subset G^IM = K^IMA^IN^I$, donc $a_1^{-1}ka_2 \in K^IM$. Par unicité du facteur $\ov{A^+}$ dans la décomposition de Cartan $G = K \ov{A^+} K$, on en déduit que $a_1=a_2=a$. Et l'élément $k$ appartient à $K^IM \cap aK^IMa^{-1} = (K^I \cap aK^Ia^{-1})M$.

\bigskip

Pour montrer la réciproque, soient $I \subsetneq \Delta$ une partie propre de $\Delta$, $a \in \ov{A^{I,+}}$ et $k_1,k_2 \in K$ tels que $k_1^{-1}k_2 \in (K^I \cap aK^Ia^{-1})M$. Alors $a^{-1}k_1^{-1}k_2a \in K^IM \subset D^I$, donc
$$ \hspace{5.4cm} k_1a D^I a^{-1} k_1^{-1} = k_2a D^I a^{-1} k_2^{-1}. \hspace{5.4cm} \qed $$

\epthm

\bigskip
\textbf{Caractérisation des groupes limites}
\bigskip

Nous commençons par retrouver la classification des sous-groupes moyennables maximaux de $G$ effectuée par Moore~\cite[Theorem~3.2, p.~133]{moore_amenable}. Notons que l'hypothèse technique de connexité isotropique n'est plus nécessaire.

\bthm[Moore] \label{thm:moyennable} Les sous-groupes moyennables fermés maximaux de $G$ sont les stabilisateurs des points de $\ov{X}^\Ch$. Ils forment donc $2^r$ classes de conjugaison (où $r=\rang (G)=\Card \Delta$), dont des représentants sont les groupes $N_G(D^I)=K^IMA_IN_I$, pour $I$ une partie de $\Delta$. \ethm

Pour démontrer ceci, nous allons nous servir du résultat suivant. Rappelons qu'un \df{plat} d'un espace métrique géodésique est un sous-espace fermé, complet, totalement géodésique et isométrique à un espace euclidien $\R^d$, où $d \geq 1$.

\bthm \label{thm:moyennable_fixe} Soit $X$ une variété de Hadamard (par exemple un espace symétrique de type non compact) et $L$ un groupe d'isométries moyennable de $X$. Alors au moins l'une des deux assertions suivantes est vraie.
\begin{enumerate}
\item Le groupe $L$ fixe un point de la compactification géométrique $\ov{X}^g$ de $X$.
\item L'espace $X$ contient un plat $L$-invariant.
\end{enumerate}
(Voir~\cite[Theorem~2, p.~506]{burger_schroeder}, ou \cite[Theorem, p.~184]{adams_ballmann}.) \hfill \qed \ethm

\bcor \label{cor:moyennable_fixe} Si $L$ est un sous-groupe moyennable fermé de $G$ non compact, alors il existe une partie propre $I$ de $\Delta$ telle qu'un conjugué de $L$ soit inclus dans $P^I=K^IMAN$. \ecor

\bp D'après le théorème~\ref{thm:moyennable_fixe}, ou bien $L$ fixe un point $\xi$ de la compactification géo\-métrique $\ov{X}^g$ de $X$, ou bien $L$ stabilise un plat $F$ de $X$.

Si $L$ fixe le point $\xi$, alors le groupe $L$ est inclus dans le stabilisateur $G_\xi$ de $\xi$ dans $G$, qui est égal à un conjugué de $P^I$, où $I$ est une partie de $\Delta$ (propre car $L$ n'est pas compact).

Si $L$ stabilise le plat $F$, alors quitte à conjuguer on peut supposer qu'il s'agit du plat $F = A_I \cdot x_0$, où $I$ est une partie de $\Delta$ (propre car $L$ n'est pas compact). Ainsi $L$ est inclus dans le stabilisateur $K^IMA_I \subset K^IMAN$ du plat $F$.
\ep

Donnons maintenant la preuve du théorème~\ref{thm:moyennable}. 

\bp Soit $L$ un sous-groupe moyennable fermé de $G$, non compact. D'après le corollaire~\ref{cor:moyennable_fixe}, considérons $I$ une partie propre minimale de $\Delta$ telle que $L \subset K^IMAN$.

\bigskip 

Faisons agir $L$ sur le sous-espace $X^I = G^I \cdot x_0$ de $X$ : il s'identifie à l'espace symétrique de type non compact associé à $G^I$. Puisque le sous-groupe $A_IN_I$ est distingué dans $K^IMAN$, le quotient $LA_IN_I / A_IN_I $ s'identifie à un sous-groupe moyennable fermé $L^I$ de $K^IMAN / A_IN_I = K^IMA^IN^I = G^IM$. Puisque le groupe $M \subset K$ normalise $G^I$, le groupe $G^IM$ agit par isométries sur le sous-espace symétrique de type non compact $X^I = G^I \cdot x_0 \subset X$. 

\bigskip

Rappelons qu'on peut choisir dans $\gg^I$ pour sous-algèbre de Cartan $\aa^I$, pour système de racines $\Sigma|_{\aa^I}$ et pour base $I|_{\aa^I}$. Supposons que $L^I$ ne soit pas compact, auquel cas d'après le corollaire~\ref{cor:moyennable_fixe} il existe une partie propre $J$ de $I$ telle que $L^I \subset K^JMA^IN^I$, quitte à conjuguer par un élément de $G^I$. Cela implique que $L$ est inclus dans $K^JMAN$, ce qui contredit la minimalité de $I$. Ainsi $L^I$ est compact, donc quitte à conjuguer par un élément de $G^I$ on peut supposer que $L^I$ est inclus dans $K^IM$. Ainsi le groupe $L$ est inclus dans $K^IMA_IN_I$.

\bigskip

Par ailleurs le groupe $K^IMA_IN_I$, extension d'un groupe compact par un groupe réso\-luble, est moyennable. Ceci achève donc la preuve.
\ep

Le théorème suivant retrouve la classification des sous-groupes distaux maximaux de $G$ (voir~\cite[Theorem~9.27, p.~144]{guivarch}), qui permet de caractériser intrinsèquement les éléments de la compactification de Chabauty.

\bthm \label{thm:distal} Les éléments de $\ov{X}^\Ch$ sont les sous-groupes distaux maximaux de $G$. Ils forment donc $2^r$ classes de conjugaison (où $r=\rang (G)=\Card \Delta$), dont des représentants sont les groupes $D^I=K^IMN_I$, pour $I$ une partie de $\Delta$. \ethm

Pour démontrer ce résultat, nous utiliserons sur le théorème suivant.

\bthm[Conze, Guivarc'h] \label{thm:distal_moy} Si $V$ est un espace vectoriel réel et si $G$ est un sous-groupe fermé distal de $\GL(V)$, alors $G$ est moyennable. \ethm

\bp
D'après le théorème~\cite[Théorème~1, p.~1083]{conze_guivarch}, on sait qu'il existe un drapeau
$$V_0 = \{0\} \subset V_1 \subset \ldots \subset V_n=V$$
stable par $G$, et tel que pour tout $i \in [[0,n-1]]$ le groupe $G$ opère de façon isométrique sur $V_{i+1}/V_i$ (pour un produit scalaire convenable).

Or le groupe des isométries linéaires d'un espace vectoriel euclidien est compact, donc le sous-groupe distingué $N$ de $G$ constitué des éléments qui agissent trivialement sur chacun des quotients $V_{i+1}/V_i$ est cocompact. Or, dans une base adaptée au drapeau, le groupe $N$ est représenté par des matrices triangulaires supérieures avec des $1$ sur la diagonale, donc est nilpotent, puis moyennable. Ainsi le groupe $G$, extension de deux groupes moyennables, est moyennable.
\ep

Donnons maintenant la preuve du théorème~\ref{thm:distal}.

\bp D'après la proposition~\ref{pro:distal}, on sait que les éléments de $\ov{X}^\Ch$ sont des sous-groupes distaux.

\bigskip

Réciproquement, soit $L$ un sous-groupe distal maximal de $G$. Puisque l'ensemble des éléments distaux de $G$ est fermé dans $G$, le sous-groupe $L$ est fermé par maximalité, donc est un sous-groupe de Lie de $G$. Les éléments de $\Ad L$ agissant sur $\gg$ ont des valeurs propres de module $1$, donc ceci reste vrai lorsqu'on restreint leur action à l'algèbre de Lie de $L$. Ainsi le groupe de Lie $L$ est distal.

\bigskip

D'après le théorème~\ref{thm:distal_moy}, on sait que le groupe $L$ est moyennable. Puis d'après le théorème~\ref{thm:moyennable}, on sait quitte à conjuguer qu'il existe une partie $I$ de $\Delta$ telle que $L \subset K^IMA_IN_I$. Fixons $l \in L$, notons sa décomposition $l=kan$, où $k \in K^IM$, $a \in A_I$ et $n \in N_I$. Puisque $K^IM$ centralise $A_I$, normalise $N_I$ et que $A_I$ normalise $N_I$, on en déduit que pour tout entier $p \in \N$ nous avons $l^p =  k^p a^p n_p$, où $n_p \in N_I$. Puisque le groupe $K^IM$ est compact et que $L$ est distal, on en déduit que le spectre de la suite $(\Ad a^pn_p)_{p \in \N}$ doit être borné. Or ce spectre est le même que le spectre de $(\Ad a^p)_{p \in \N}$, donc il ne peut être borné que si $a=1$. En conclusion, on en déduit que $L$ est inclus dans le sous-groupe distal $K^IMN_I = D^I$, donc $L=D^I$ par maximalité.
\ep

Remarquons que tout sous-groupe fermé distal de $G$ est à croissance polynomiale, mais que la réciproque n'est pas vraie. Par exemple le sous-groupe abélien $A$ est à croissance polynomiale, mais n'est pas distal en tant que sous-groupe de $G$ : il est distal en tant que groupe de Lie.

\section{Lien avec la compactification polyédrale}

\bigskip
\textbf{Compactification polyédrale d'un espace vectoriel}
\bigskip

Deux références pour la compactification polyédrale sont~\cite[\S~I.2, p.~21]{landvogt} et \cite{remy}.

\bigskip

Soit $V$ un espace vectoriel réel de dimension finie $r$. Soit $\Sigma$ un ensemble fini non vide de formes linéaires non nulles sur $V$, stable par $\alpha \mapsto -\alpha$, tel qu'on obtienne une partition de $V \bs \{0\}$ indexée par les partitions $\Sigma=\Sigma_0 \sqcup \Sigma_+ \sqcup \Sigma_-$ de $\Sigma$ en trois parties de $\Sigma$, en cônes simpliciaux ouverts (facettes) définis par 
$$F_{\Sigma_0,\Sigma_+,\Sigma_-} = \{x \in V\bs \{0\} : \forall \alpha \in \Sigma_0, \alpha(x)=0,\; \forall \alpha \in \Sigma_+, \alpha(x)>0 \mbox{ et } \forall \alpha \in \Sigma_-, \alpha(x)<0 \}.$$ 

Notons $\mathcal{F}$ l'ensemble des facettes, et $\mathcal{F}_0$ l'ensemble des chambres (facettes de dimension maximale). Si $F$ est une facette de cette décomposition, on note $\langle F \rangle$ le sous-espace vectoriel de $V$ engendré par $F$. On définit la \df{compactification polyédrale} de $V$ relativement à $\Sigma$ comme l'ensemble union disjointe
$$\ov{V}^p = \bigsqcup_{F \in \mathcal{F}} V/ \langle F \rangle ,$$
muni de la topologie que nous allons définir.

Si $F_0 \in \mathcal{F}_0$ est une chambre, on appelle \df{coin} associé à $F_0$ le sous-ensemble de $\ov{V}^p$ :
$$C_{F_0} = \bigsqcup_{F \in \mathcal{F},\; F \subset \ov{F_0}} V/ \langle F \rangle ,$$
où $\ov{F_0}$ désigne l'adhérence de $F_0$ dans $V$. Remarquons que le coin $C_{F_0}$ contient une copie de $V$, pour la facette $\{0\}$.

Soient $\alpha_1,...,\alpha_r \in \Sigma$ tels que $F_0 = \{x \in V : \forall i \in [[1,r]], \alpha_i(x)>0\}$. Ceci permet de considérer l'application :
\beq \label{eqn:topo_coin} \beta : C_{F_0} & \ra & ]-\infty,+\infty]^r \\
x \in V/ \langle F \rangle & \mapsto & (\beta(x)_i)_{i \in [[1,r]]} \eeq
$$ \mbox{où } \beta(x)_i = \left\{ \begin{array}{cc} \alpha_i(x) & \mbox{ si } {\alpha_i}|_F = 0 \\
+\infty & \mbox{ si } {\alpha_i}|_F > 0 \end{array} \right. .$$

L'application $\beta$ est une bijection : munissons $C_{F_0} \subset \ov{V}^p$ de la topologie qui fait de la bijection $\beta$ un homéomorphisme de $C_{F_0}$ sur $]-\infty,+\infty]^r$, muni de la topologie usuelle.

\bigskip

On munit $\ov{V}^p$ de la topologie faible définie par la famille $(C_{F_0})_{F \in \mathcal{F}_0}$~: une partie $U \subset \ov{V}^p$ est ouverte si et seulement si, pour toute chambre $F_0 \in \mathcal{F}_0$, la partie $U \cap C_{F_0}$ est ouverte dans $C_{F_0}$~: ceci fait de $\ov{V}^p$ un espace compact. Si l'on note $W$ le groupe des automorphismes linéaires de $V$ qui préservent l'ensemble $\mathcal{F}$ des facettes, alors l'espace $\ov{V}^p$ est une $W$-compactification de $V$ (voir~\cite[\S~I.2, p.~21]{landvogt}).

\bigskip
\textbf{Compactification polyédrale d'un espace symétrique}
\bigskip

Soit $X$ un espace symétrique de type non compact, et $x_0 \in X$ un point base. Soit $G$ un groupe de Lie connexe agissant continûment isométriquement sur $X$, se surjectant sur $\Isom_0(X)$ avec noyau fini. Soit $K =G_{x_0}$ le stabilisateur du point $x_0$. Soit $\gg = \kk \oplus \aa \oplus \nn$ une décomposition d'Iwasawa de l'algèbre de Lie $\gg$ de $G$. Le système de racines $\Sigma$ associé à une sous-algèbre de Lie $\aa$ de $\gg$ abélienne diagonalisable maximale donne une décomposition de $\aa \bs \{0\}$ en cônes simpliciaux. 

\bigskip

D'après ce qui précède, on peut définir la compactification polyédrale $\ov{\aa}^p$ de $\aa$. Si on choisit une chambre de Weyl (ouverte) positive $\aa^+$, notons alors $\ov{\aa^+}^p$ l'adhérence de $\aa^+$ dans $\ov{\aa}^p$. Notons $\Delta$ la base du système de racines $\Sigma$ associée à la chambre de Weyl positive $\aa^+$. Les facettes de la chambre de Weyl $\aa^+$ sont alors exactement les $\aa_I^+$, où $I$ parcourt les parties propres de $\Delta$. Ainsi les facettes ajoutées pour obtenir la compactification polyédrale de $\ov{\aa^+}$ sont les $\aa / \aa_I$, où $I$ parcourt les parties propres de $\Delta$.

\bigskip

Remarquons que l'application :
\beq K \times \ov{\aa^+} & \ra & X \\
(k,H) & \mapsto & k\exp(H) \cdot x_0 \eeq
est continue et surjective. Montrons qu'elle est propre : soit $(k_n\exp(H_n) \cdot x_0)_{n \in \N}$ une suite dans l'image de cette application convergeant vers $k\exp(H) \cdot x_0$. Par compacité de $K$, on peut supposer quitte à extraire que la suite $(k_n)_{n \in \N}$ converge vers $k' \in K$. Alors la suite $(\exp(H_n) \cdot x_0)_{n \in \N}$ converge vers $k'^{-1}k\exp(H) \cdot x_0$ dans $X$. Ainsi, il existe une suite $(k'_n)_{n \in \N}$ d'éléments de $K$ telle que la suite $(\exp(H_n)k'_n)_{n \in \N}$ converge vers $k'^{-1}k\exp(H)$ dans $G$. Par continuité de la composante dans $\ov{A^+}$ dans le décomposition de Cartan $G=K\ov{A^+}K$, on en déduit que la suite $(\exp(H_n))_{n \in \N}$ converge vers $\exp(H)$ dans $\ov{A^+}$, et donc que la suite $(H_n)_{n \in \N}$ converge vers $H$ dans $\ov{\aa^+}$. En conclusion, la suite $(k_n,H_n)_{n \in \N}$ converge vers $(k',H)$ dans $K \times \ov{\aa^+}$. On a donc montré que l'application $K \times \ov{\aa^+} \ra X$ est propre.

Ainsi cette application, par passage au quotient, fournit un homéomorphisme dont on notera l'inverse
$$ \psi : X \ra (K \times \ov{\aa^+}) / \sim ,$$
où $(k,H) \sim (k',H')$ si $k\exp(H) \cdot x_0 = k'\exp(H') \cdot x_0$.

\bigskip

Ceci suggère d'étendre la relation d'équivalence à $G \times \ov{\aa^+}^p$, en posant, si $H \in \aa^+$ et $H' \in \aa^+$~:
$$(g,H) \sim (g',H') \Longleftrightarrow g\exp(H)K\exp(-H)g^{-1} = g'\exp(H')K\exp(-H')g'^{-1}.$$
Et, si $I$ et $I'$ sont des parties propres de $\Delta$ et si $H + \aa_I \in \aa / \aa_I$ et $H' + \aa_{I'} \in \aa / \aa_{I'}$ :
$$(g,H + \aa_I) \sim (g',H' + \aa_{I'}) \Longleftrightarrow g\exp(H)D^I\exp(-H)g^{-1} = g'\exp(H')D^{I'}\exp(-H')g'^{-1}.$$
Cette définition a bien un sens car, pour toute partie propre $I$ de $\Delta$, le groupe $A_I$ normalise $D^I$. En effet, le groupe $A$ normalise le groupe $N_I$ et centralise le groupe $M$. De plus, le groupe $A_I$ centralise le groupe $G^I$, donc en particulier centralise le groupe $K^I$. 

\bigskip

Si $I=\Delta$, alors $\aa_\Delta=\{0\}$ donc on peut identifier $H + \aa_\Delta$ et $H \in \aa^+$. De plus, si l'on note $D^\Delta = K$, ceci permet d'écrire la relation d'équivalence sous la même forme pour tous les éléments de $G \times \ov{\aa^+}^p$ : 
$$(g,H + \aa_I) \sim (g',H' + \aa_{I'}) \Longleftrightarrow g\exp(H)D^I\exp(-H)g^{-1} = g'\exp(H')D^{I'}\exp(-H')g'^{-1} ,$$
où $I$ et $I'$ sont des parties quelconques de $\Delta$.

\bigskip

\bpro \label{pro:interversion} L'application $f:\ov{\aa^+}^p \ra \Ch(G)$ définie par $H+\aa_I \mapsto \exp(H) D^I \exp(-H)$ est continue. \epro

\bp  Les espaces $\ov{\aa^+}^p$ et $\Ch(G)$ étant métrisables, nous allons montrer que $f$ est séquen\-tiellement continue : soit $(H_n+\aa_{I_n})_{n \in \N}$ une suite de $\ov{\aa^+}^p$ convergeant vers $H_\infty+\aa_J$, où $H_n \in \aa^{I_n}$ et $H_\infty \in \aa^J$. Le nombre de facettes étant fini, on peut supposer quitte à extraire que la partie $I_n \subset \Delta$ est constante, égale à $I$. D'après la topologie sur $\ov{\aa^+}^p$, nous avons nécessairement $J \subset I$. Or, par continuité de l'exponentielle et de l'action de $G$ par conjugaison sur $\Ch(G)$, l'application $f$ restreinte à la facette $\aa^J$ est continue : ainsi, on peut supposer que $J$ est strictement inclus dans $I$. Quitte à translater par $H_\infty$, on peut supposer que $H_\infty + \aa_J=\aa_J$.

\bigskip

Soit $U$ un voisinage fermé de $D^J$ dans $\Ch(G)$ : montrons que $\exp(H_n) D^I \exp(-H_n)$ appartient à $U$ si $n$ est assez grand, ce qui montrera la continuité de $f$. D'après la proposition~\ref{pro:di}, il existe des constantes $\eps>0$ et $N>0$ telles que pour tout $H$ appartenant à
$$V_{\eps,N} = \{H \in \aa^+ \,:\, \forall \alpha \in \Delta \bs J, \alpha(H) \geq N \mbox{ et } \alpha \in J, |\alpha(H)| \leq \eps\},$$
nous ayons $\exp(H)K\exp(-H) \in U$.

\bigskip

Puisqu'on peut choisir les $H_n$ modulo $\aa_I$, on peut supposer que pour tout $\alpha \in \Delta \bs I$ nous ayons $\alpha(H_n)=0$. Puisque la suite $(H_n+\aa_I)_{n \in \N}$ converge vers $\aa_J$ dans $\ov{\aa^+}^p$, on en déduit que pour tout $\alpha \in I \bs J$ nous avons $\alpha(H_n) \ral{n \ra +\infty} +\infty$, et que pour tout $\alpha \in J$ nous avons $\alpha(H_n) \ral{n \ra +\infty} 0$. Ainsi il existe un entier $n_0 \in \N$ tel que pour tout $n \geq n_0$ nous ayons $H_n \in V_{\eps,N}$.

\bigskip

Soit $(H'_m)_{m \in \N}$ une suite de $\aa_I^+$ telle que pour tout $\alpha \in \Delta \bs I$ nous ayons $\alpha(H'_m) \ral{n \ra +\infty} +\infty$. Ainsi pour tout $m \in \N$ et tout $n \geq n_0$ nous avons les deux propriétés suivantes :
\beq \forall \alpha \in \Delta \bs J, & \alpha(H_n+H'_m) \geq \alpha(H_n) \geq N \\
\forall \alpha \in J, & |\alpha(H_n+H'_m)| = |\alpha(H_n)|\leq \eps. \eeq
C'est-à-dire que $H_n+H'_m \in V_{\eps,N}$, donc $\exp(H_n+H'_m)K\exp(-H_n-H'_m) \in U$. Or lorsque l'on fait tendre $m$ vers $+\infty$, on en déduit d'après la proposition~\ref{pro:di} que, pour tout $n \geq n_0$ nous avons
$$ \liml_{m \ra +\infty} \exp(H_n+H'_m)K\exp(-H_n-H'_m) = \exp(H_n)D^I\exp(-H_n) \in U.$$
Ceci montre la continuité de $f$.\ep

Remarquons que dans~\cite[Remark~7.34, p.~115]{guivarch}, la continuité de l'action de $G$ sur la compactification polyédrale est démontrée via la compactification de Martin.

\bigskip

Définissons la \df{compactification polyédrale} de $X$ par
$$\ov{X}^p = (G \times \ov{\aa^+}^p) / \sim ,$$
où l'espace $\ov{X}^p$ est muni de la topologie quotient, et où l'inclusion $X \ra \ov{X}^p$ est la composition de l'homéomorphisme $ \psi : X \ra (K \times \ov{\aa^+}) / \sim$ et de l'inclusion $(K \times \ov{\aa^+}) /\!\sim \; \ra (G \times \ov{\aa^+}^p) /\!\sim$.

L'espace $\ov{X}^p$ est naturellement muni d'une $G$-action à gauche, quotient de l'action de $G$ à gauche sur le premier facteur de $G \times \ov{\aa^+}^p$.

\bigskip

Si $G^0$ désigne la composante neutre du groupe des isométries de $X$, alors le morphisme $\pi : G \ra G^0$ fourni par l'action de $G$ sur $X$ est, par hypothèse, surjectif et de noyau fini. Les compactifications polyédrales de $X$ obtenues pour $G$ et $G^0$ sont naturellement isomorphes, par l'homéomorphisme $\pi$-équivariant
\beq (G \times \ov{\aa^+}^p) /\!\sim & \ra & (G^0 \times \ov{\aa^+}^p) /\!\sim \\
\left[(g,H)\right] & \mapsto & \left[(\pi(g),H)\right] .\eeq

\bpro \label{pro:Gcompactification} L'espace $\ov{X}^p$ est une $G$-compactification de $X$. \epro

\bp Puisque l'application $f$ est continue, l'application $G \times \ov{\aa^+}^p \ra \Ch(G)$ qui à $(g,H+\aa_I)$ associe $g\exp(H)D^I\exp(-H)g^{-1}$ est continue et passe au quotient. Ainsi, l'espace quotient $\ov{X}^p$ est séparé. La projection $K \times \ov{\aa^+}^p \ra \ov{X}^p$ est surjective (tout élément de $G \times \ov{\aa^+}^p$ est équivalent à un élément de $K \times \ov{\aa^+}^p$) et continue (par passage au quotient de l'inclusion $K \times \ov{\aa^+}^p \ra G \times \ov{\aa^+}^p$), or l'espace $\ov{X}^p$ est séparé et l'espace $K \times \ov{\aa^+}^p$ est compact, on conclut donc que $\ov{X}^p$ est compact.

\bigskip

Montrons que l'espace $X$, identifié par l'homéomorphisme $\psi:X \ra (K \times \ov{\aa^+}) / \! \sim$ à un sous-espace de $\ov{X}^p$, est dense dans $\ov{X}^p$. Soit $(g,H+\aa_I)$ un point de $G \times \ov{\aa^+}^p$ tel que son image $x$ dans $\ov{X}^p$ n'appartienne pas à $\psi(X)$, c'est-à-dire tel que $I$ soit une partie propre de $\Delta$. Tout voisinage de $(g,H+\aa_I)$ dans $G \times \ov{\aa^+}^p$ intersecte $G \times \aa^+$, car par définition $\ov{\aa^+}^p$ est l'adhérence de $\aa^+$ dans $\ov{\aa}^p$. Ainsi, tout voisinage de $x$ dans $\ov{X}^p$ intersecte $\psi(X)$. Par conséquent, l'espace $\ov{X}^p$ est bien une compactification de l'espace symétrique $X$.

\bigskip

L'homéomorphisme $X \ra G \times \ov{\aa^+} / \! \sim$ est $G$-équivariant, donc l'espace $\ov{X}^p$ est une $G$-compactification de $X$.
\ep

\bpro \label{pro:poly_chabauty} Les $G$-compactifications de Chabauty $\ov{X}^\Ch$ et polyédrale $\ov{X}^p$ de l'espace symétrique de type non compact $X$ sont $G$-isomorphes, via le $G$-homéomorphisme $\varphi : \ov{X}^p \ra \ov{X}^\Ch$ qui à la classe de $(g,H+\aa_I)$ associe $g\exp(H)D^I\exp(-H)g^{-1}$. \epro

\bp D'après la définition de la relation d'équivalence $\sim$, l'application $\varphi$ est bien définie et injective. D'après le début de la preuve de la proposition~\ref{pro:Gcompactification}, l'application $\varphi$ est continue. D'après le théorème~\ref{thm:existence_ecriture}, l'application $\varphi$ est surjective. Or l'espace $\ov{X}^p$ est compact et $\ov{X}^\Ch$ est compact donc séparé : ainsi $\varphi$ est un homéomorphisme, $G$-équivariant. \ep

\bpro \label{pro:poly_cell} La compactification polyédrale $\ov{X}^p$ est $G$-isomorphe à la compactification cellulaire duale $X \cup \Delta^*(X)$, telle qu'elle est définie dans~\cite[Definition~3.40, p.~45]{guivarch}. \epro

\bp Nous allons vérifier que la compactification polyédrale $\ov{X}^p$ satisfait les hypothèses du théorème~\cite[Theorem~3.38, p.~43]{guivarch}, qui caractérise la compactification cellulaire duale.

\bigskip

Soit $(k_n,H_n)_{n \in \N}$ une suite fondamentale au sens de~\cite[Definition~3.35, p.~41]{guivarch}, c'est-à-dire qu'il existe une partie propre $I$ de $\Delta$ telle que la suite $(k_n)_{n \in \N}$ de $K$ converge vers $k$, et si $H_n = H_{n,I}+H_n^I$ est la décomposition de $H_n \in \aa$ selon $\aa = \aa_I \oplus \aa^I$, alors la suite $(H_{n,I})_{n \in \N}$ tend vers $+\infty$ dans $\aa_I^+$ et la suite $(H_n^I)_{n \in \N}$ converge vers $H \in \ov{\aa^{I,+}}$. Alors, d'après les propositions~\ref{pro:di} et \ref{pro:poly_chabauty}, la suite $(k_n\exp H_n)_{n \in \N}$ converge, dans la compactification polyédrale, vers la classe d'équivalence de $(k,H+\aa_I)$.

\bigskip

Soient $(k_n,H_n)_{n \in \N}$ et $(k'_n,H'_n)_{n \in \N}$ deux suites fondamentales. Leurs limites formelles au sens de~\cite[Definition~3.35, p.~41]{guivarch} sont $(\Ad k(\aa_I^+),k \exp(H) \cdot x_0)$ et $(\Ad k'(\aa_{I'}^+),k' \exp(H') \cdot x_0)$ respectivement. Leurs limites dans la compactification polyédrale sont les classes d'équivalence de $(k,H+\aa_I)$ et $(k',H'+\aa_{I'})$ respectivement. Notons $a = \exp(H)$ et $a' = \exp(H')$.

Supposons que les limites formelles soient égales. Alors, d'après la proposition~\cite[Proposition~3.4, p.~25]{guivarch}, nous avons $I=I'$ et $k^{-1}k$ appartient au normalisateur $Z_K(\aa_I) = K^IM$ de $\aa_I$ dans $K$. De plus $ka \cdot x_0 = k'a' \cdot x_0 \in M$, donc $a^{-1}k^{-1}k'a' \in K$. Or $a$ et $a'$ appartiennent à $\ov{A^+}$, donc par unicité de la décomposition de Cartan $G=K\ov{A^+}K$, nous avons $a=a'$. Par ailleurs $a^{-1}k^{-1}k'a \in A^IK^IMA^I \cap K \subset G^IM \cap K = K^IM$. Par conséquent $k^{-1}k' \in aK^IMa^{-1} \cap K^IM$. Donc, d'après la réciproque du théorème~\ref{thm:unicite_ecriture}, on en déduit que
$$ka D^I a^{-1}k^{-1} = k'a' D^{I'} a'^{-1} k'^{-1}.$$
Donc les limites dans la compactification polyédrale des classes d'équivalence de $(k,H+\aa_I)$ et $(k',H+\aa_I)$ sont égales.

Réciproquement, supposons que les limites dans la compactification polyédrale soient égales. Alors, d'après le théorème~\ref{thm:unicite_ecriture}, on en déduit que $I=I'$, $a=a'$ et $k^{-1}k' \in (K^I \cap aK^Ia^{-1})M$. Alors $k^{-1}k'$ appartient au stabilisateur $aKa^{-1}$ du point $a \cdot x_0$, donc
$$k' \exp(H') \cdot x_0 = k'a \cdot x_0 = ka \cdot x_0 = k \exp(H) \cdot x_0 .$$
De plus $k^{-1}k' \in K^IM = Z_K(\aa_I)$, donc d'après la proposition~\cite[Proposition~3.4, p.~25]{guivarch}, nous avons $\Ad k(\aa_I^+) = \Ad k'(\aa_I^+)$. Ainsi, les limites formelles sont égales.

\bigskip

Ainsi, la compactification polyédrale $\ov{X}^p$ satisfait les hypothèses du théorème~\cite[Theorem~3.38, p.~43]{guivarch}, donc est isomorphe à la compactification cellulaire duale. Ces compactifications étant des $G$-compactifications de $X$ (voir~\cite[Theorem~3.43, p.~46]{guivarch}), ce sont des compactifications $G$-isomorphes de $X$. \ep

\bthm Les compactifications de Chabauty, polyédrale, cellulaire duale (\cite[Definition~3.40, p.~45]{guivarch}), de Satake-Furstenberg maximale (\cite[Proposition~4.53, p.~73]{guivarch}) et de Martin (\cite[Chapter~VII, pp.~103--115]{guivarch}) sont $G$-isomorphes. \ethm

\bp D'après les résultats~\cite[Theorem~4.43, p.66]{guivarch}, \cite[Proposition~4.53, p.~73]{guivarch} et \cite[Theorem~7.33, p.~115]{guivarch}, les compactifications cellulaire duale, de Satake-Furstenberg maximale et de Martin sont $G$-isomorphes. D'après les propositions~\ref{pro:poly_chabauty} et \ref{pro:poly_cell}, celles-ci sont également $G$-isomorphes aux compactifications de Chabauty et polyédrale. \ep



\begin{thebibliography}{99}

\bibitem{abels_affine}
Abels, H., \emph{Distal affine transformation groups}, J.~Reine~Angew.~Math. \textbf{299/300} (1978), 294--300.

\bibitem{abels_lie}
Abels, H., \emph{Distal automorphism groups of Lie groups}, J.~Reine~Angew.~Math. \textbf{329} (1981), 82--87.

\bibitem{adams_ballmann}
Adams, S., et W.~Ballmann, \emph{Amenable isometry groups of Hadamard spaces}, Math.~Ann. \textbf{312} (1998), 183--195.

\bibitem{borel_ji}
Borel, A., et L.~Ji, ``Compactifications of symmetric spaces'', J.~Diff.~Geom. \textbf{75} (2007), 1--56.

\bibitem{bourbaki}
Bourbaki, N., \og \'El\'ements de math\'ematique. Int\'egration \fg, Masson-Dunod, 1959.

\bibitem{bridson_haefliger}
Bridson, M.~R. et A.~Haefliger, ``Metric spaces of non-positive curvature'', Grund. math.~Wiss., Springer Verlag, 1999.

\bibitem{heisenberg}
Bridson, M.~R., P.~de~la~Harpe, et V.~Kleptsyn, \emph{The Chabauty space of closed subgroups of the three-dimensional Heisenberg group}, arXiv:0711.3736, à paraître dans Pacific J.~Math. (2007).

\bibitem{burger_schroeder}
Burger, M. et V.~Schroeder, \emph{Amenable groups and stabilizers of measures on the boundary of a Hadamard manifold}, Math.~Ann. \textbf{276} (1987), 505--514.

\bibitem{CEG}
Canary, R.D., D.B.A.~Epstein et P.L.~Green, \emph{Notes on notes of Thurston}, dans ``Analytical and geometric aspects of hyperbolic space'', D.B.A.~Epstein ed., Lond.~Math.~Soc., Lect.~Notes~Series {\bf 111} (1987), 3--92.

\bibitem{chabauty}
Chabauty, C., \emph{Limite d'ensembles et g\'eom\'etrie des nombres}, Bull.~Soc.~Math.~France \textbf{78} (1950), 143--151.

\bibitem{conze_guivarch}
Conze, J.-P. et Y.~Guivarc'h, \emph{Remarques sur la distalité dans les espaces vectoriels}, CRAS Paris \textbf{278} (1974), 1083--1086.

\bibitem{cdp}
Courtois, G., F.~Dal'bo et F.~Paulin, \og Sur la dynamique des groupes de matrices et applications arithmétiques \fg, Journées mathématiques X-UPS, Les \'{E}ditions de l'\'{E}cole Polytechnique,  2007.

\bibitem{eberlein}
Eberlein, P.B., ``Geometry of nonpositively curved manifolds'', Chicago~Lect.~Math.~, Univ.~Chicago Press, 1996.

\bibitem{ghl}
Gallot, S., D.~Hulin et J.~Lafontaine, ``Riemannian geometry'', Springer Verlag, 1987.

\bibitem{guivarch}
Guivarc'h, Y., L.~Ji et J.C.~Taylor, ``Compactifications of symmetric spaces", Progr.~Math. {\bf 156}, 1998.

\bibitem{guivarch_eigen}
Guivarc'h, Y., L.~Ji et J.C.~Taylor, \emph{Compactifications of symmetric spaces and positive eigenfunctions of the Laplacian}, CRM Proc.~Lect.~Notes \textbf{28} (2001), 69--116.

\bibitem{m2}
Haettel, T., \emph{Compactification de Chabauty des espaces symétriques de type non compact}, mémoire de seconde année de mastère, Univ.~Orsay (2008).

\bibitem{RxZ}
Haettel, T., \emph{L'espace des sous-groupes fermés de $\R \times \Z$}, prépublication ENS Ulm, arXiv:0905.3290 (2008), à paraître dans Algebraic and Geometric Topology.

\bibitem{harpe_chabauty}
de~la~Harpe, P., \emph{Spaces of closed subgroups of locally compact groups}, Tripode \textbf{14}, arXiv : 0807.2030v2 (2008).

\bibitem{helgason}
Helgason, S., ``Differential geometry, Lie groups, and symmetric spaces'', Grad.~Stud.~Math. {\bf 34}, Amer.~Math.~Soc., 1978.

\bibitem{kloeckner}
Kloeckner, B., \emph{The space of closed subgroups of $\R^n$ is stratified and simply connected}, Journal of Topology 2 \textbf{3} (2009), 570--588.

\bibitem{landvogt}
Landvogt, E., ``A compactification of the Bruhat-Tits building'', Lect.~Notes~Math. {\bf 1619}, Springer Verlag, 1995.

\bibitem{mneimne_testard}
Mneimn\'e, R. et F.~Testard, \og Introduction \`a la th\'eorie des groupes de Lie classiques \fg, Hermann, 1986.

\bibitem{moore_amenable}
Moore, C.C., \emph{Amenable subgroups of semi-simple groups and proximal flows}, Isr.~J.~Math. \textbf{34} (1979), 121--138.

\bibitem{mostow}
Mostow, G.D., ``Strong rigidity of locally symmetric spaces'', Annals of Math. Studies \textbf{78}, Princeton Univ.~Press, 1973.

\bibitem{onishchik_vinberg}
Onishchik, A.L. et E.B.~Vinberg, ``Lie groups and algebraic groups'', Springer Verlag, 1990.

\bibitem{pourezza}
Pourezza, I., et J.~Hubbard, \emph{The space of closed subgroups of $\R^2$}, Topology \textbf{18} (1979), 143--146.

\bibitem{ratcliffe}
Ratcliffe, J.G., ``Foundations of hyperbolic manifolds'', Grad.~Texts~Math. {\bf 149}, Springer Verlag, 2006.

\bibitem{remy}
Rémy, B., \emph{Immeuble à l'infini et combinatoire des groupes, Compactification polyédrique des plats}, http://math.univ-lyon1.fr/\tild remy/CpctPol.pdf (1999).

\bibitem{satake}
Satake, I., \emph{On representations and compactifications of symmetric Riemannian spaces}, Ann.~Math. \textbf{71} (1960), 77--110.

\bibitem{warner}
Warner, G.,``Harmonic analysis on semi-simple Lie groups I'', Grund.~math.~Wiss. {\bf 188}, Springer Verlag, 1972.

\end{thebibliography}
\end{document}